\documentclass[12pt,a4paper]{article}

\usepackage{amsfonts,amsmath,amsthm,amssymb,amsopn}
\usepackage{graphicx}
\usepackage{epstopdf}
\usepackage{algorithmic}
\usepackage{enumitem}

\usepackage{bm,color}
\usepackage{xcolor}
\usepackage{stmaryrd}
\usepackage[T1]{fontenc}
\usepackage{float, caption, subcaption}
\usepackage{booktabs}
\usepackage{multirow}
\usepackage{siunitx}
\usepackage{pifont}
\usepackage{hyperref}
\usepackage{cleveref}
\usepackage{bbold}

\newcommand*{\I}{\mathbb{i}}

\newcommand\bu{\bm{u}}

\newcommand\dd{\mathrm{d}}

\newcommand\abs[1]{\lvert #1 \rvert}

\newcommand\xr{i+\frac12}
\newcommand\xl{i-\frac12}
\newcommand\yr{j+\frac12}
\newcommand\yl{j-\frac12}

\theoremstyle{plain}

\theoremstyle{definition}
\newtheorem{example}{Example}[section]
\newtheorem{proposition}{Proposition}[section]
\newtheorem{remark}{Remark}[section]

\crefname{equation}{Equation}{Equations}
\crefname{figure}{Figure}{Figures}
\crefname{table}{Table}{Tables}
\crefname{example}{Example}{Examples}
\crefname{section}{Section}{Sections}

\renewcommand{\title}{An asymptotic-preserving active flux scheme for the hyperbolic heat equation in the diffusive scaling}

\newcommand{\authorOne}{Junming Duan\footnote{Institute of Mathematics, University of W\"urzburg, Emil-Fischer-Stra\ss e 40, 97074 W\"urzburg, Germany, junming.duan@uni-wuerzburg.de}}
\newcommand{\authorTwo}{Wasilij Barsukow\footnote{Institut de Math\'ematiques de Bordeaux (IMB), CNRS UMR 5251, University of Bordeaux, 33405 Talence, France, wasilij.barsukow@math.u-bordeaux.fr}}
\newcommand{\authorThree}{Christian Klingenberg\footnote{Institute of Mathematics, University of W\"urzburg, Emil-Fischer-Stra\ss e 40, 97074 W\"urzburg, Germany, christian.klingenberg@uni-wuerzburg.de}}

\begin{document}

\begin{center} \Large
\title

\vspace{1cm}

\date{}
\normalsize

\authorOne, \authorTwo, \authorThree
\end{center}

\begin{abstract}

The Active Flux (AF) method is a compact, high-order finite volume scheme that enhances flexibility by introducing point values at cell interfaces as additional degrees of freedom alongside cell averages.
The method of lines is employed here for temporal discretization.
A common approach for updating point values relies on the Jacobian Splitting (JS) method, which incorporates upwinding.
A key advantage of the AF method over standard finite volume schemes is its structure-preserving property, motivating the investigation of its asymptotic-preserving (AP) behavior in the diffusive scaling.
We show that the JS-based AF method without any modification is AP for solving the hyperbolic heat equation, in the sense that the limit scheme is a discretization of the limit heat equation. We use formal asymptotic analysis, discrete Fourier analysis, and numerical experiments to illustrate our findings.




\end{abstract}

\tableofcontents

\section{Introduction}\label{sec:introduction}

Particle systems in physics, such as rarefied gases and neutron transport, can be modeled at different scales.
At a microscopic level, particles move and collide individually, and the system is often described by kinetic models,
which provide a statistical description using probability distributions.
At a macroscopic level, the mean free path, or the average distance between collisions, is small compared to a macroscopic length,
and the particles behave like fluids.
In this case, macroscopic models such as fluid or diffusion equations offer accurate and computationally efficient approximations.
The passage from microscopic to macroscopic descriptions can be formally derived via asymptotic analysis.

Numerically, multiscale problems where the mean free path varies across several orders of magnitude pose significant challenges.
Traditional kinetic solvers require grid resolutions smaller than the mean free path, making computations prohibitively expensive.
To address this, asymptotic preserving (AP) schemes have been developed to seamlessly transition between kinetic and macroscopic models.
Their stability and convergence are independent of the mean free path,
and allow coarse grids, which are crucial for applications such as inertial confinement fusion.
The AP schemes were first proposed in \cite{Larsen_1987_Asymptotic_JoCP,Larsen_1989_Asymptotic_JoCP} for the steady-state solutions to neutron transport equation in the diffusive regime and have been extended to a broad range of kinetic models.
Some AP finite volume schemes were constructed based on careful modification of numerical flux \cite{Jin_1996_Numerical_JoCP} and/or upwinding treatment of source term \cite{Gosse_2002_asymptotic_CRM}.
Different reformulations were proposed for designing AP schemes,
e.g., even-odd decomposition by rewriting linear transport equations as parity equations \cite{Klar_1998_asymptotic_SJoNA, Jin_2000_Uniformly_SJoNA}.
Another popular way is micro-macro decomposition by decomposing the density distribution function into the local Maxwellian plus the deviation,
e.g., for the radiative heat transfer equations \cite{Klar_2001_Numerical_MMaMiAS},
linear kinetic equations \cite{Lemou_2008_new_SJoSC},
nonlinear Boltzmann equation \cite{Bennoune_2008_Uniformly_JoCP}, and so on.
High-order AP discontinuous Galerkin (DG) schemes were developed using micro-macro decomposition and several different numerical fluxes \cite{Jang_2015_High_JoCP,Peng_2020_Stability_JoCP}.
For a more comprehensive review, we refer readers to \cite{Jin_2022_Asymptotic_AN} and references therein.

The active flux (AF) method is a new compact finite volume method \cite{Eymann_2011_Active_InCollection, Eymann_2011_Active_InProceedings,Eymann_2013_Multidimensional_InCollection,Roe_2017_Is_JoSC},
inspired by \cite{VanLeer_1977_Towards_JoCP}.
It simultaneously evolves cell averages and additional degrees of freedom, chosen as point values located at cell interfaces.
Thanks to these point values, the AF method does not need Riemann solvers (unlike Godunov methods), as the numerical solution is continuous across the cell interface.
The AF methods can be roughly divided into two classes based on the evolution of the point value.
The original ones evolve the cell average through Simpson's rule for flux quadrature in time,
and employ exact or approximate evolution operators to evolve the point values and to obtain the numerical solutions at the flux quadrature points. Examples are exact evolution operators for linear equations \cite{Barsukow_2019_Active_JoSC,Fan_2015_Investigations_InCollection,Eymann_2013_Multidimensional_InCollection, VanLeer_1977_Towards_JoCP},
$p$-system \cite{Fan_2017__AcousticComponent},
and approximate evolution operators for Burgers' equation \cite{Eymann_2011_Active_InCollection,Eymann_2011_Active_InProceedings,Roe_2017_Is_JoSC,Barsukow_2021_active_JoSC},
the compressible Euler equations in one spatial dimension \cite{Eymann_2011_Active_InCollection,Helzel_2019_New_JoSC,Barsukow_2021_active_JoSC},
multidimensional Euler equations \cite{Fan_2017__AcousticComponent},
and hyperbolic balance laws \cite{Barsukow_2021_Active_SJoSC, Barsukow_2023_Well_CoAMaC}, etc.
The method of bicharacteristics was used for the derivation of truly multidimensional approximative evolution operators \cite{Chudzik_2024_Active_JoSC}.
The other so-called generalized, or semi-discrete AF methods adopt a method of lines,
where the evolution of the cell average and point value is written in semi-discrete form and integrated in time by using Runge-Kutta methods. Examples of this approach are \cite{Abgrall_2023_Combination_CoAMaC, Abgrall_2023_Extensions_EMMaNA, Abgrall_2025_Semi_JoSC, Abgrall_2024_Bound} based on Jacobian splitting (JS) and \cite{Duan_2025_Active_SJoSC, Duan_2025_Active} based on flux vector splitting.

The AF method is superior to standard finite volume methods  due to its structure-preserving property.
For multi-dimensional acoustic equations, it preserves the vorticity and stationary states  \cite{Barsukow_2019_Active_JoSC},
and for acoustics with gravity, it is naturally well-balanced \cite{Barsukow_2021_Active_SJoSC}.
These encouraging results lead us to the central question of this paper: Does the AF method also possess the AP property in the diffusive scaling?
To start, we examine the hyperbolic heat equation, the lowest-order angular discretization of the transport equation.
This equation is also known as the $P_1$ model, the first-order formulation of the telegraph equation, or the damped wave equation \cite{Buet_2012_Design_NM}.
Since it is a linear system, a fully-implicit time discretization is adopted.
We present formal asymptotic analysis and discrete Fourier analysis for the JS-based AF scheme for solving the hyperbolic heat equation,
and show that it is AP.
Several 1D and 2D numerical results verify our theoretical findings.
The key feature of such JS-based AF method is that it is automatically AP without any modification.

The remainder of this paper is structured as follows.
\Cref{sec:hyperbolic_heat} introduces the hyperbolic heat equation and its limit heat equation by using formal asymptotic analysis.
\Cref{sec:1d_af_scheme} gives the 1D AF scheme based on the JS for the point value update, and \Cref{sec:1d_af_formal} derives the limit scheme via formal asymptotic analysis.
As we observe order degradation for the JS-based AF scheme in the limit, an alternative point value update is discussed in \Cref{sec:alternating}.
A discrete Fourier analysis is adopted to study the convergence order of the scheme and its limit in \Cref{sec:fourier}. 
\Cref{sec:2d_af_scheme} discusses the 2D case.
Numerical tests are conducted in \Cref{sec:results} to experimentally demonstrate the accuracy and AP property.
\Cref{sec:conclusion} concludes the paper with final remarks.

\section{Hyperbolic heat equation and its limit in the diffusive scaling}\label{sec:hyperbolic_heat}

The hyperbolic heat equation in the diffusive scaling reads
\begin{equation}\label{eq:hyperbolic_heat}
	\left\{~
	\begin{aligned}
		&p_t + \frac{1}{\epsilon}\nabla\cdot\bu = 0, \\
		&\bu_t + \frac{1}{\epsilon}\nabla p + \frac{\sigma}{\epsilon^2}\bu = \bm{0},
	\end{aligned}
	\right.
\end{equation}
where $0<\epsilon<1$ and $\sigma>0$.

To study the asymptotic limit of \eqref{eq:hyperbolic_heat} as $\epsilon\rightarrow 0$,
we assume sufficient regularity of the solutions and perform a formal analysis by expanding them as a power series in $\epsilon$,
\begin{align}\label{eq:power_series}
	p &= p^{(0)} + \epsilon p^{(1)} + \epsilon^2 p^{(2)} + \dots,&
	\bu &= \bu^{(0)} + \epsilon \bu^{(1)} + \epsilon^2 \bu^{(2)} + \dots.
\end{align}
The asymptotic behavior is determined by inserting \eqref{eq:power_series} into \eqref{eq:hyperbolic_heat}.
Matching the terms of the same order, we have
\begin{align*}
	\epsilon^{-2}: ~&\bu^{(0)} = \bm{0}, \\
	\epsilon^{-1}: ~&\nabla p^{(0)} + \sigma \bu^{(1)} = \bm{0}, \\
	\epsilon^{0}: ~&p^{(0)}_t + \nabla\cdot\bu^{(1)} = 0.
\end{align*}
Thus, the leading order solutions satisfy
\begin{equation*}
	p^{(0)}_t = \nabla\cdot \left(\frac{1}{\sigma}\nabla p^{(0)}\right), \quad \bu^{(0)} = \bm{0},
\end{equation*}
i.e., $p^{(0)}$ is the solution to the limit heat equation.
\section{1D active flux scheme}\label{sec:1d_af_scheme}
This section presents the 1D semi-discrete AF methods for the 1D hyperbolic heat equation,
\begin{equation}\label{eq:1d_hyperbolic_heat}
	\left\{~
	\begin{aligned}
		&p_t + \frac{1}{\epsilon}u_x = 0, \\
		&u_t + \frac{1}{\epsilon}p_x + \frac{\sigma}{\epsilon^2}u = 0.
	\end{aligned}
	\right.
\end{equation}
As the advection and source term are both stiff as $\epsilon\ll 1$, and the system is linear, the fully-discrete methods are obtained using the 3rd-order four-stage Diagonally Implicit Runge-Kutta (DIRK) method from \cite{Kennedy_2016_Diagonally, Kennedy_2019_Diagonally_ANM},
which is stiffly accurate.

\subsection{Update of cell average}
Assume that a 1D computational domain is divided into $N$ uniform cells
$I_i = [x_{\xl}, x_{\xr}]$ with the cell size $\Delta x=x_{\xr}-x_{\xl}$ and cell centers $x_i = (x_{\xl} + x_{\xr})/2$. 
The degrees of freedom consist of the cell averages and point values at the cell interfaces,
\begin{align*}
	&\bar{p}_i(t) = \frac{1}{\Delta x} \int_{I_i} p_h(x,t) \dd x,\quad
	\bar{u}_i(t) = \frac{1}{\Delta x} \int_{I_i} u_h(x,t) \dd x, \\
	&p_{\xr}(t) = p_h(x_{\xr}, t),\quad u_{\xr}(t) = u_h(x_{\xr}, t),
\end{align*}
where $p_h(x,t)$ and $u_h(x,t)$ are numerical approximations.
As for finite volume methods, the update of the cell average is found upon integrating \eqref{eq:1d_hyperbolic_heat} over $I_i$,
\begin{align}
	&\frac{\dd}{\dd t}\bar{p}_i = - \frac{1}{\epsilon\Delta x}\left(u_{\xr} - u_{\xl}\right), \label{eq:1d_av_p} \\
	&\frac{\dd}{\dd t}\bar{u}_i = - \frac{1}{\epsilon\Delta x}\left(p_{\xr} - p_{\xl}\right) - \frac{\sigma}{\epsilon^2}\bar{u}_i. \label{eq:1d_av_u}
\end{align}
Note that the flux at the cell interface is available directly, without the need for a Riemann solver.

\subsection{Update of point value}
For the point value evolution, the Jacobian splitting (JS) \cite{Abgrall_2023_Extensions_EMMaNA} is adopted.
The Jacobian matrix of \eqref{eq:1d_hyperbolic_heat} can be diagonalized,
\begin{equation*}
	\bm{J} = 
	\frac{1}{\epsilon}
	\begin{pmatrix}
		0 & 1 \\ 1 & 0 \\
	\end{pmatrix}
	= \begin{pmatrix}
		1 & -1 \\ 1 & 1 \\
	\end{pmatrix}
	\begin{pmatrix}
		\frac{1}{\epsilon} & 0 \\ 0 & -\frac{1}{\epsilon} \\
	\end{pmatrix}
	\begin{pmatrix}
		\frac{1}{2} & \frac{1}{2} \\ -\frac{1}{2} & \frac{1}{2} \\
	\end{pmatrix}
	=: \bm{R}\bm{\Lambda}\bm{R}^{-1},
\end{equation*}
and split into positive and negative parts as
\begin{equation*}
	\bm{J}^+ = \bm{R}\bm{\Lambda}^{+}\bm{R}^{-1}
	= \frac{1}{\epsilon}\begin{pmatrix}
		\frac12 \phantom{\Big|}& \frac12 \\ \frac12\phantom{\Big|} & \frac12 \\
	\end{pmatrix},\quad
	\bm{J}^- = \bm{R}\bm{\Lambda}^{-}\bm{R}^{-1}
	= \frac{1}{\epsilon}\begin{pmatrix}
		-\frac12 & \phantom{\Big|}\frac12 \\ \phantom{\Big|}\frac12 & -\frac12 \\
	\end{pmatrix},
\end{equation*}
such that $\bm{\Lambda}^\pm = \frac12\left(\bm{\Lambda} \pm \abs{\bm{\Lambda}}\right)$.
Then the point value update can be written as
\begin{align*}
	\frac{\dd}{\dd t}\begin{pmatrix}
		p_{\xr} \\ u_{\xr} \\
	\end{pmatrix}
	= -\left[\bm{J}^{+}D^{+}_{\xr}\begin{pmatrix}
		p \\ u \\
	\end{pmatrix}
	+ \bm{J}^{-}D^{-}_{\xr}\begin{pmatrix}
		p \\ u \\
	\end{pmatrix}\right]
	- \frac{\sigma}{\epsilon^2}\begin{pmatrix}
		0 \\ u_{\xr} \\
	\end{pmatrix},
\end{align*}
where $D^{+}_{\xr}$ and $D^{-}_{\xr}$ are component-wise upwind finite difference operators.
These are obtained by first reconstructing a parabola $u_\texttt{para}(x)$ in the cell $I_i$ satisfying
\begin{equation*}
	\frac{1}{\Delta x}\int_{I_i} u_\texttt{para}(x)\dd x = \bar{u}_i,\quad
	u_\texttt{para}(x_{i\pm\frac12}) = u_{i\pm\frac12},
\end{equation*}
whose derivative at $x_{\xr}$ yields a high-order accurate upwind finite difference
\begin{equation*}
	D_{\xr}^+(u) = u_\texttt{para}'(x_{\xr}) = \frac{1}{\Delta x}\left(2u_{\xl} - 6\bar{u}_{i} + 4u_{\xr}\right).
\end{equation*}
For negative eigenvalues, the reconstruction in cell $I_{i+1}$ is differentiated at $x_{\xr}$ instead, yielding
\begin{equation*}
	D_{\xr}^-(u) = \frac{1}{\Delta x}\left(- 4u_{\xr} + 6\bar{u}_{i+1} - 2u_{i+\frac32}\right).
\end{equation*}
Therefore, the point value update can be rewritten as
\begin{align}
	&\frac{\dd}{\dd t} p_{\xr} = - \frac{1}{2\epsilon}\Big[(D_{\xr}^+ - D_{\xr}^-)p + (D_{\xr}^+ + D_{\xr}^-)u\Big], \label{eq:1d_pnt_p} \\
	&\frac{\dd}{\dd t} u_{\xr} = - \frac{1}{2\epsilon}\Big[(D_{\xr}^+ + D_{\xr}^-)p + (D_{\xr}^+ - D_{\xr}^-)u\Big] - \frac{\sigma}{\epsilon^2} u_{i+\frac12}. \label{eq:1d_pnt_u}
\end{align}

\section{Formal asymptotic analysis}\label{sec:1d_af_formal}

Expand sufficiently smooth numerical solutions as a power series in $\epsilon$,
\begin{align*}
	\bar{p}_i &= \bar{p}_i^{(0)} + \epsilon \bar{p}_i^{(1)} + \epsilon^2 \bar{p}_i^{(2)} + \dots,&
	\bar{u}_i &= \bar{u}_i^{(0)} + \epsilon \bar{u}_i^{(1)} + \epsilon^2 \bar{u}_i^{(2)} + \dots, \\
	p_{\xr} &= p_{\xr}^{(0)} + \epsilon p_{\xr}^{(1)} + \epsilon^2 p_{\xr}^{(2)} + \dots,&
	u_{\xr} &= u_{\xr}^{(0)} + \epsilon u_{\xr}^{(1)} + \epsilon^2 u_{\xr}^{(2)} + \dots,
\end{align*}
and substitute them into the schemes \eqref{eq:1d_av_p}-\eqref{eq:1d_pnt_u}.
After equating the coefficients of equal powers of $\epsilon$, one has
\begin{align*}
	\epsilon^{-2}:\quad &\bar{u}_i^{(0)} = u_{\xr}^{(0)} = 0, \\
	\epsilon^{-1}:\quad &\bar{u}_i^{(1)} = -\frac{1}{\sigma\Delta x}\left(p^{(0)}_{\xr} - p^{(0)}_{\xl}\right), \\
	&u_{\xr}^{(1)} = -\frac{1}{2\sigma}\Big[(D_{\xr}^+ + D_{\xr}^-)p^{(0)} + (D_{\xr}^+ - D_{\xr}^-)u^{(0)}\Big] \\
	&\quad\quad = -\frac{1}{2\sigma}\left(D_{\xr}^+ + D_{\xr}^-\right)p^{(0)}, \\
	\epsilon^{0}:\quad &\frac{\dd}{\dd t} \bar{p}_{i}^{(0)} = -\frac{1}{\Delta x}\left(u^{(1)}_{\xr} - u^{(1)}_{\xl}\right), \\ 
	&\frac{\dd}{\dd t} p_{\xr}^{(0)} = -\frac{1}{2}\Big[(D_{\xr}^+ - D_{\xr}^-)p^{(1)} + (D_{\xr}^+ + D_{\xr}^-)u^{(1)}\Big].
\end{align*}
The second equation at order $\epsilon^{-1}$ implies 
\begin{align*}
	-u^{(1)}_{\xr} =&\ \frac{1}{2\sigma}(D_{\xr}^+ + D_{\xr}^-)p^{(0)}
	=\ \frac{1}{\sigma\Delta x}\left( p_{i-\frac12}^{(0)} - 3\bar{p}_{i}^{(0)} + 3\bar{p}_{i+1}^{(0)} - p_{i+\frac32}^{(0)}\right) \\
	=&\ \frac{1}{\sigma}\left( p^{(0)}_{x}\Big|_{x_{\xr}} \!\!- \, \frac{\Delta x^2}{12}p^{(0)}_{xxx}\Big|_{x_{\xr}}
	+ \mathcal{O}(\Delta x^4)\right),
\end{align*}
thus the first equation at order $\epsilon^0$ can be simplified as
\begin{align*}
	\frac{\dd}{\dd t} \bar{p}_{i}^{(0)} 
	=&\ \frac{1}{\sigma\Delta x} p^{(0)}_{x}\Big|_{x_{\xl}}^{x_{\xr}}
	- \frac{\Delta x}{12\sigma}
	p^{(0)}_{xxx}\Big|_{x_{\xl}}^{x_{\xr}}
	+ \mathcal{O}(\Delta x^3) \\
	=&\ \frac{1}{\sigma\Delta x} p^{(0)}_{x}\Big|_{x_{\xl}}^{x_{\xr}}
	+ \mathcal{O}(\Delta x^2),
\end{align*}
which is a finite volume approximation of $p^{(0)}_t = \frac{1}{\sigma}p^{(0)}_{xx}$ with the truncation error $\mathcal{O}(\Delta x^2)$.
Note that the last equality uses the Lipschitz continuity of $p^{(0)}_{xxx}$ if $p^{(0)}$ is sufficiently smooth.
The second equation of the order $\epsilon^0$ can be simplified as
\begin{align*}
	\frac{\dd}{\dd t} p_{\xr}^{(0)} =& -\frac{1}{2}\Big[(D_{\xr}^+ - D_{\xr}^-)p^{(1)} + (D_{\xr}^+ + D_{\xr}^-)u^{(1)}\Big] \\
	=& -\frac{1}{\sigma\Delta x}\Big[ u^{(1)}_{\xl} - 3\bar{u}^{(1)}_{i} + 3\bar{u}^{(1)}_{i+1} - u^{(1)}_{i+\frac32} \Big] +\frac{1}{2}\Big[(D_{\xr}^- - D_{\xr}^+)p^{(1)}\Big]  \\
	=&\ \frac{1}{\sigma\Delta x}\Big[ \frac12(D_{\xl}^+ + D_{\xl}^-)p^{(0)} - 3\left( p^{(0)}_{\xr} - p^{(0)}_{\xl} \right) + 3\left( p^{(0)}_{i+\frac32} - p^{(0)}_{\xr} \right) \\
	& - \frac12(D_{i+\frac32}^+ + D_{i+\frac32}^-)p^{(0)} \Big] +\frac{1}{2}\Big[(D_{\xr}^- - D_{\xr}^+)p^{(1)}\Big] \\
	=&\ \frac{1}{\sigma\Delta x^2}\Big(p_{i-\frac32}^{(0)} - 3\bar{p}_{i-1}^{(0)} + 3p_{i-\frac12}^{(0)} + 3\bar{p}_{i}^{(0)} - 8p_{i+\frac12}^{(0)} + 3\bar{p}_{i+1}^{(0)} + 3p_{i+\frac32}^{(0)} - 3\bar{p}^{(0)}_{i+2} + p_{i+\frac52}^{(0)} \Big) \\
	&+ \frac{1}{\Delta x} \left(- p_{i-\frac12}^{(1)}+3\bar{p}_i^{(1)} - 4p_{i+\frac12}^{(1)} +3\bar{p}_{i+1}^{(1)}-p_{i+\frac32}^{(1)}\right) \\
	=&\ \frac{1}{\sigma}\left( p^{(0)}_{xx}\Big|_{x_{i+\frac12}}
	+ \frac{1}{12}p_{xxxx}^{(0)}\Big|_{x_{i+\frac12}}\Delta x^2\right)
	-\frac{1}{30}p_{xxxx}^{(1)}\Big|_{x_{i+\frac12}}\Delta x^3 + \mathcal{O}(\Delta x^4),
\end{align*}
whose right-hand side recovers an approximation of $\frac{1}{\sigma}p^{(0)}_{xx}$ with the truncation error $\mathcal{O}(\Delta x^2)$.

\begin{proposition}
	Assuming sufficient regularity of the solutions, as $\epsilon\rightarrow 0$, the leading order solutions of the 1D JS-based AF schemes \eqref{eq:1d_av_p}-\eqref{eq:1d_pnt_u} satisfy
	\begin{align}
		\bar{u}_i^{(0)} =&\ u_{\xr}^{(0)}= 0, \nonumber \\
		\frac{\dd}{\dd t} \bar{p}_{i}^{(0)} =&\ \frac{1}{2\sigma\Delta x}\left[ (D_{\xr}^+ + D_{\xr}^-)p^{(0)} - (D_{\xl}^+ + D_{\xl}^-)p^{(0)} \right], \label{eq:1d_limit_av_p0} \\
		\frac{\dd}{\dd t} p_{\xr}^{(0)} =&\ \frac{1}{\sigma\Delta x}\Big[ \frac12(D_{\xl}^+ + D_{\xl}^-)p^{(0)} - 3\left( p^{(0)}_{\xr} - p^{(0)}_{\xl} \right) + 3\left( p^{(0)}_{i+\frac32} - p^{(0)}_{\xr} \right) \nonumber \\
		& - \frac12(D_{i+\frac32}^+ + D_{i+\frac32}^-)p^{(0)} \Big]
		+ \mathcal{O}(\Delta x^3), \label{eq:1d_limit_pnt_p0}
	\end{align}
	where \eqref{eq:1d_limit_av_p0}-\eqref{eq:1d_limit_pnt_p0} are approximations of $p^{(0)}_t = \frac{1}{\sigma}p^{(0)}_{xx}$ with the truncation errors $\mathcal{O}(\Delta x^2)$.
	To be specific,
	\begin{equation*}
		\frac{\dd}{\dd t} \bar{p}_{i}^{(0)} = \frac{1}{\sigma\Delta x} p^{(0)}_{x}\Big|_{x_{\xl}}^{x_{\xr}}
		+ \mathcal{O}(\Delta x^2),\quad
		\frac{\dd}{\dd t} p_{\xr}^{(0)} =
		\frac{1}{\sigma} p^{(0)}_{xx}\Big|_{x_{i+\frac12}} + \mathcal{O}(\Delta x^2).
	\end{equation*}
\end{proposition}

The truncation error may not fully explain the convergence order of the AF scheme (or other compact schemes, such as Discontinuous Galerkin, see \cite{Zhang_2003_MMMAS_AnalysisThree}).
Indeed, if one considers the most classical AF scheme for the linear advection equation,
the truncation error apparently is 2nd order while the overall accuracy is actually 3rd order (see e.g. \cite{Abgrall_2023_Extensions_EMMaNA,Hasan_2023_JoCaAM_CentralCompact,Zeng_2014_C&F_HighorderHybrid}).
To examine more precisely the accuracy of the limit schemes,
a Fourier analysis is performed in Section \ref{sec:fourier}, after discussing an alternative point value update next.

\section{An alternative point value update}\label{sec:alternating}

The following is another way to update the point values:
\begin{align*}
	&\frac{\dd}{\dd t} p_{\xr} = - \frac{1}{\epsilon} D_{\xr}^+ u, \\
	&\frac{\dd}{\dd t} u_{\xr} = - \frac{1}{\epsilon} D_{\xr}^- p - \frac{\sigma}{\epsilon^2} u_{i+\frac12}.
\end{align*}
It is inspired by the ``alternating flux'' in \cite{Jang_2015_High_JoCP}.
This choice of upwinding is very different from how stable methods are usually achieved for hyperbolic systems of PDEs: Characteristics are not considered and instead, a left bias is chosen on one variable and a right bias on another. We do not find this prescription to yield a stable method in two space dimensions. The reason for discussing it are its theoretical properties that will be of importance in the subsequent analysis of Section \ref{sec:fourier}.

Using the formal asymptotic analysis, we have
\begin{align*}
	\epsilon^{-2}:\quad &\bar{u}_i^{(0)} = u_{\xr}^{(0)} = 0, \\
	\epsilon^{-1}:\quad &\bar{u}_i^{(1)} = -\frac{1}{\sigma\Delta x}\left(p^{(0)}_{\xr} - p^{(0)}_{\xl}\right), \\
	&u_{\xr}^{(1)} = -\frac{1}{\sigma} D_{\xr}^- p^{(0)}, \\
	\epsilon^{0}:\quad &\frac{\dd}{\dd t} \bar{p}_{i}^{(0)} = -\frac{1}{\Delta x}\left(u^{(1)}_{\xr} - u^{(1)}_{\xl}\right), \\ 
	&\frac{\dd}{\dd t} p_{\xr}^{(0)} = - D_{\xr}^+ u^{(1)}.
\end{align*}
Thus the limit scheme is
\begin{align*}
    \frac{\dd}{\dd t} \bar{p}_{i}^{(0)} &= \frac{1}{\sigma\Delta x^2}\left(D_{\xr}^- p^{(0)} - D_{\xl}^- p^{(0)}\right), \\
    \frac{\dd}{\dd t} p_{\xr}^{(0)} &= - \frac{1}{\Delta x} \left( 2u_{\xl}^{(1)} - 6\bar{u}_i^{(1)} + 4u_{\xr}^{(1)} \right) \\
    & = \frac{1}{\sigma\Delta x^2} \left[ 2D_{\xl}^- p^{(0)} - 6\left(p_{\xr}^{(0)} - p_{\xl}^{(0)}\right) + 4D_{\xr}^-p^{(0)} \right].
\end{align*}

\section{Fourier analysis}\label{sec:fourier}
\subsection{Fourier analysis of the PDEs}

We first perform a Fourier analysis for the hyperbolic heat equation \eqref{eq:1d_hyperbolic_heat}.
Note that the Fourier analysis is also used in \cite{Jin_1996_Numerical_JoCP} to study the parabolic behavior of the AP schemes.
Make the following Fourier ansatz of the solutions
\begin{equation*}
	p(x,t) = \hat{p}(t)\exp(\I \omega x), \quad u(x,t) = \hat{u}(t)\exp(\I \omega x),
\end{equation*}
and substitute them into \eqref{eq:1d_hyperbolic_heat},
then we have the following evolution equation
\begin{equation}
	\frac{\dd}{\dd t}
	\begin{bmatrix}
		\hat{p} \\ \hat{u} \\
	\end{bmatrix}
	= \underbrace{-\frac{1}{\epsilon}\begin{bmatrix}
		0 & \I\omega \\
		\I\omega & \epsilon^{-1}\sigma \\
	\end{bmatrix}}_{=:\mathcal{E}}
	\begin{bmatrix}
		\hat{p} \\ \hat{u} \\
	\end{bmatrix}. \label{eq:evoFourier}
\end{equation}

The eigenvalues of the evolution matrix $\mathcal{E}$ are
\begin{equation*}
 \lambda_{1,2} = \frac{-\sigma \pm s}{2 \epsilon^2}
\end{equation*}
with $s := \sqrt{\sigma^2 - 4 \epsilon^2 \omega^2}$, which leads to
\begin{align}
	\lambda_1 &= -\frac{\omega^2}{\sigma} - \frac{\omega^4 \epsilon^2}{\sigma^3} + \mathcal{O}(\epsilon^4),&
	\lambda_2 &= -\frac{\sigma}{\epsilon^2} + \frac{\omega^2}{\sigma}  + \frac{\omega^4 \epsilon^2}{\sigma^3} + \mathcal{O}(\epsilon^4), \label{eq:physicaleval}
\end{align}

Given the initial condition
\begin{equation}\label{eq:1d_fourier_init}
	p_0 = \hat{p}_0\exp(\I\omega x),\quad u_0 = \hat{u}_0\exp(\I\omega x),
\end{equation}
one can decompose $[\hat{p}_0,\hat{u}_0]^\top = V_1 + V_2$ in the eigenvectors of $\mathcal E$ as follows
\begin{equation*}
V_1 = \left[ \hat p_0 + \mathcal O(\epsilon),~ \mathcal O(\epsilon) \right]^\top,\quad
V_2 = \left[ \mathcal O(\epsilon), ~ \hat u_0 + \mathcal O(\epsilon)\right]^\top.
\end{equation*}

The solution of \eqref{eq:evoFourier} is
\begin{equation*}
	[\hat{p}, \hat{u}]^\top	
	= \exp(\lambda_1 t) V_1 + \exp(\lambda_2 t) V_2.
\end{equation*}
In the limit $\epsilon\rightarrow 0$, the first eigenvalue $\lambda_1$ tends to $-\frac{1}{\sigma}\omega^2$,
which corresponds to the solution of the limit heat equation for $p$, as can be seen from $V_1$.
The second eigenvalue diverges as $\epsilon \to 0$ and corresponds to the initial layer, occurring if the initial data are not well-prepared, i.e., if $\hat u_0 \neq 0$ as is obvious from the nature of $V_2$.
The corresponding evolution $\exp(\lambda_2 t) = \exp\left( - \frac{\sigma t}{\epsilon^2}\right)$ decays rapidly with time if $\epsilon$ is small. Seen as a function of $\epsilon$, this term has a Taylor series uniformly zero, and hence escapes the formal analysis presented in Section \ref{sec:hyperbolic_heat}.
For later reference, we state here the eigenvectors 
\begin{align}
V_1 &= \frac{1}{2s} \Big[ \hat p_0 (s+\sigma) - 2 \I \epsilon \omega \hat u_0, ~\hat u_0 (s-\sigma) - 2 \I \epsilon \omega \hat p_0 \Big ],
\label{eq:eigenvectorpde1} \\
V_2 &=  \frac{1}{2s}\Big[ \hat p_0 (s-\sigma) + 2 \I \epsilon \omega \hat u_0, ~\hat u_0 (s+\sigma) + 2 \I \epsilon \omega \hat p_0 \Big ],
\label{eq:eigenvectorpde2}
\end{align}
where $s - \sigma = \mathcal O(\epsilon^2)$.

\begin{remark}
    Note that $s$ may be complex.
    In fact, for $\epsilon=0$ it is real corresponding to the parabolic limit, while for $\sigma=0$, one has $s = \pm 2\I \epsilon\omega$ and $\lambda_{1,2} = \pm \I\omega/\epsilon$, i.e., the imaginary eigenvalues expected in the hyperbolic regime.
\end{remark}


\subsection{Fourier analysis of the numerical scheme}

Inspired by \cite{Zhang_2003_MMMAS_AnalysisThree, Guo_2013_JoCP_SuperconvergenceDiscontinuous}, we now perform a similar\footnote{The calculations are performed using Mathematica.} Fourier analysis for the AF scheme by taking the following Fourier ansatz of the numerical solutions
\begin{align}
	\bar{p}_i(t) &= \hat{p}_1(t)\exp(\I \omega x_i),\quad & p_{\xr}(t) &= \hat{p}_2(t)\exp(\I \omega x_i), \label{eq:discretefourieransatz1}\\
	\bar{u}_i(t) &= \hat{u}_1(t)\exp(\I \omega x_i),\quad & u_{\xr}(t) &= \hat{u}_2(t)\exp(\I \omega x_i), \label{eq:discretefourieransatz2}
\end{align}
where $\omega$ is the wave number\footnote{It is of no relevance whether we associate the point values with $x_i$ or $x_{i+\frac12}$, as the latter would amount to an extra overall factor $\exp(\I \omega \Delta x/2)$ and hence only a redefinition of $\hat p_2$ or $\hat u_2$.}.
Let also $t_x=\exp(\I\omega\Delta x)$.

The coefficient $\hat w_0 = [\hat p_1(0), \hat u_1(0), \hat p_2(0), \hat u_2(0)]^\top$ of the initial data $\hat w_0\exp(\I\omega x_i)$ for the numerical solution in the cell $I_i=[x_{i-\frac12}, x_{i+\frac12}]$ is obtained based on \eqref{eq:1d_fourier_init},
\begin{align}
	\hat{w}_0  & = \exp(-\I \omega x_i) \left[\frac{1}{\Delta x}\int_{I_i} p_0(x) \dd x, ~ \frac{1}{\Delta x}\int_{I_i} u_0(x) \dd x, ~ p_0(x_{\xr}), ~ u_0(x_{\xr})\right]^\top \nonumber\\
    &= \left[\hat p_0 \frac{2 \sin\left(\frac{\omega \Delta x }{2}\right)}{\omega \Delta x} , ~\hat u_0 \frac{2 \sin\left(\frac{\omega \Delta x }{2}\right)}{\omega \Delta x} , ~\hat p_0 \exp\left( \I \omega \tfrac{\Delta x}{2} \right), ~\hat u_0 \exp\left( \I \omega \tfrac{\Delta x}{2} \right) \right]^\top. \label{eq:projfouriermodeontoDof}
\end{align}

\subsubsection{Point value update using Jacobian splitting}

Substitute the ansatz \eqref{eq:discretefourieransatz1}--\eqref{eq:discretefourieransatz2} into the AF schemes \eqref{eq:1d_av_p}-\eqref{eq:1d_pnt_u}:
\begin{equation*}
	\frac{\dd}{\dd t} [\hat{p}_1, \hat{u}_1, \hat{p}_2, \hat{u}_2]^\top
	= \bm{G} [\hat{p}_1, \hat{u}_1, \hat{p}_2, \hat{u}_2]^\top, 
\end{equation*}
where the evolution matrix is
\begin{equation*}
	\bm{G}
	= \frac{1}{\epsilon\Delta x}\begin{bmatrix}
		0 & 0 & 0 & t_x^{-1}-1 \\
		0 & -\epsilon^{-1}\sigma\Delta x & t_x^{-1}-1 & 0 \\
		3(1+t_x) & 3(1-t_x) & -(t_x^{-1} + 4 + t_x) & t_x-t_x^{-1} \\
		 3(1-t_x) & 3(1+t_x) & t_x-t_x^{-1} & -\epsilon^{-1}\sigma\Delta x - (t_x^{-1} + 4 + t_x) \\
	\end{bmatrix}.
\end{equation*}
It is difficult to calculate the eigenvalues of $\bm{G}$,
and we thus construct them as a power series in $\Delta x$
\begin{equation*}
    \tilde{\lambda} =  \Delta x^{-1}\tilde{\lambda}^{(-1)} + \tilde{\lambda}^{(0)} + \Delta x \tilde{\lambda}^{(1)} + \Delta x^2 \tilde{\lambda}^{(2)} + \Delta x^3 \tilde{\lambda}^{(3)} + \text{h.o.t.},
\end{equation*}
and solve the characteristic equation
\begin{equation*}
	\det (\tilde{\lambda} \bm{I}_4 - \bm{G}) = 0
\end{equation*}
by matching the coefficients of the same order in 
$\Delta x$.
The eigenvalues are
\begin{align*}
	\tilde{\lambda}_1 &= \lambda_1 - \frac{\Delta x^3 \omega^4}{72 \epsilon} + \mathcal{O}(\Delta x^4), \\
    \tilde{\lambda}_2 &= \lambda_2 - \frac{\Delta x^3 \omega^4}{72 \epsilon} + \mathcal{O}(\Delta x^4), \\
	\tilde{\lambda}_3 &= -\frac{6}{\epsilon\Delta x} - \frac{\tilde s + \sigma}{2 \epsilon^2} + \frac{\omega^2\Delta x }{\epsilon} - \frac{2 \Delta x^2 \omega^4}{\tilde s} - \frac{5 \Delta x^3 \omega^4}{72 \epsilon} + \mathcal O(\Delta x^4), \\
	\tilde{\lambda}_4 &= -\frac{6}{\epsilon\Delta x} + \frac{\tilde s - \sigma}{2 \epsilon^2} + \frac{\omega^2\Delta x }{\epsilon} + \frac{2 \Delta x^2 \omega^4}{\tilde s} - \frac{5 \Delta x^3 \omega^4}{72 \epsilon} + \mathcal O(\Delta x^4),
\end{align*}
where $\tilde s := \sqrt{\sigma^2 - 36 \epsilon^2 \omega^2}$.
The first and second eigenvalues are physical (i.e., converging to \eqref{eq:physicaleval} as $\Delta x \to 0$), with the other two spurious modes decaying exponentially as $\Delta x \to 0$.
Note that the second physical eigenvalue corresponds to the initial layer.
Observe also the presence of a singularity in the physical eigenvalues as $\epsilon \to 0$; we will come back to this aspect in Section \ref{ssec:discussion} after presenting the rest of the analysis.

Similarly, we calculate a set of eigenvectors $\tilde v_1, \dots, \tilde v_4$ of $\bm{G}$ by assuming that they are power series of $\Delta x$.
Then the discrete initial data $\hat{w}_0$ can be decomposed in the basis of the eigenvectors of $\bm{G}$ as $\hat{w}_0 = \hat{w}_0^{(1)} \tilde v_1 + \hat{w}_0^{(2)}\tilde v_2 + \hat{w}_0^{(3)}\tilde v_3 +\hat{w}_0^{(4)} \tilde v_4 =: \tilde V_1 + \tilde V_2 + \tilde V_3 +\tilde V_4$ with, ignoring terms $\mathcal O(\Delta x^4)$,
\begin{align}
    \tilde V_1 &= \Big[ -\frac{\omega  \epsilon  (2 p_0 \omega  \epsilon +\I u_0
   (s-\sigma ))}{s (s-\sigma )}
   + \frac{\Delta x^2 \omega ^3 \epsilon  (2 p_0 \omega  \epsilon +\I
   u_0 (s-\sigma ))}{24 s (s-\sigma )}, \nonumber \\
   &\qquad -\frac{(s-\sigma )(\I p_0 \omega  \epsilon  +\sigma  u_0)+2 u_0 \omega ^2 \epsilon ^2}{s (s-\sigma )}
   + \frac{\Delta x^2 \omega ^2 ((s-\sigma )(\I p_0 \omega  \epsilon 
   +\sigma  u_0)+2 u_0 \omega ^2
   \epsilon ^2)}{24 s (s-\sigma )}, \nonumber \\
   &\qquad -\frac{\omega  (288 p_0 \omega  \epsilon ^2+144 \I u_0
   \epsilon  (s-\sigma ))}{144 s (s-\sigma )}
   -\frac{\Delta x \omega  (-72 u_0 \omega  \epsilon 
   (s-\sigma )+144 \I p_0 \omega ^2 \epsilon ^2)}{144 s
   (s-\sigma )} \nonumber\\
   &\qquad -\frac{\Delta x^2 \omega  (-36 p_0 \omega ^3 \epsilon ^2-18
   \I u_0 \omega ^2 \epsilon  (s-\sigma ))}{144 s (s-\sigma
   )} \nonumber\\
   &\qquad -\frac{1}{144 s (s-\sigma )}\Big( \Delta x^3 \omega  (2 p_0 \omega ^3 \epsilon  (s-\sigma )-6
   \I p_0 \omega ^4 \epsilon ^2-2 \I \sigma  u_0 \omega ^2
   (s-\sigma ) \nonumber\\
   &\qquad +3 u_0 \omega ^3 \epsilon  (s-\sigma )-4 \I
   u_0 \omega ^4 \epsilon ^2) \Big), \nonumber \\
   &\qquad \frac{-288 \I p_0 \omega  \epsilon +144 (s-\sigma) u_0}{288 s}
   + \frac{\Delta x (144 p_0 \omega ^2 \epsilon +72 \I (s-\sigma)
   u_0 \omega)}{288 s} \nonumber\\
   &\qquad + \frac{\Delta x^2 (36 \I p_0 \omega ^3 \epsilon -18 (s-\sigma)
   u_0 \omega ^2)}{288 s} \nonumber\\
   &\qquad +\frac{
   \Delta x^3 (-2 \I p_0 \omega ^3 (s+\sigma )-6 p_0
   \omega ^4 \epsilon -u_0 \omega ^3 (4 \omega  \epsilon +3 \I
   s)+3 \I \sigma  u_0 \omega ^3)
   }{288s}
   \Big ]^\top, \label{eq:jsdecomp1}
\end{align}
\begin{align}
\tilde V_2 &=  \Big[\ 
\frac{\I \epsilon  \omega  (2 \I p_0 \epsilon  \omega +u_0 (s+\sigma))}{s(s+\sigma)}-\frac{\I \Delta x^2 \epsilon  \omega ^3 (2 \I p_0 \epsilon  \omega +u_0 (s+\sigma))}{24 s (s+\sigma)}, \nonumber\\
   &\qquad \frac{u_0 \sigma +\I p_0 \epsilon \omega +\frac{1}{2} u_0 (s-\sigma )}{s}-\frac{\Delta x^2 \omega ^2 (u_0 \sigma +\I p_0
   \epsilon  \omega +\frac{1}{2} u_0 (s-\sigma ))}{24 s},\nonumber\\
   &\qquad \frac{\omega  \left(-2 p_0 s \omega  \epsilon ^2+\I \sigma u_0 \epsilon  (s+\sigma )-4 \I u_0 \omega ^2 \epsilon^3\right)}{s^2 (s+\sigma )}
   + \frac{\Delta x \omega  \left(-2 \I p_0 s \omega ^2 \epsilon^2-\sigma  u_0 \omega  \epsilon  (s+\sigma )+4 u_0   \omega ^3 \epsilon ^3\right)}{2 s^2 (s+\sigma )} \nonumber\\
   &\qquad +\frac{\Delta x^2 \omega  \left(2 p_0 s \omega ^3 \epsilon ^2-\I   \sigma  u_0 \omega ^2 \epsilon  (s+\sigma )+4 \I u_0   \omega ^4 \epsilon ^3\right)}{8 s^2 (s+\sigma )} \nonumber\\
   &\qquad + \frac{1}{144 s^2 (s+\sigma )}\Big(\Delta x^3 \omega  (2 p_0 \omega ^3 \epsilon    (\sigma ^2+s \sigma +\omega  \epsilon  (-4 \omega  \epsilon +3   \I s))-2 \I \sigma ^2 u_0 \omega ^2 (s+\sigma ) \nonumber\\
   &\qquad +4 \I
   u_0 \omega ^4 \epsilon ^2 (s+2 \sigma )+3 \sigma  u_0   \omega ^3 \epsilon  (s+\sigma )-12 u_0 \omega ^5 \epsilon
   ^3)\Big), \nonumber \\
   & \qquad \frac{2 \I p_0 \omega  \epsilon + (s+\sigma) u_0}{2s}   + \frac{\Delta x (-2 p_0 \omega ^2 \epsilon +\I (s+\sigma) u_0 \omega )}{4 s} \nonumber\\
   &\qquad - \frac{\Delta x^2 (2 \I p_0 \omega ^3 \epsilon +(s+\sigma) u_0
   \omega ^2)}{16 s} \nonumber\\
   &\qquad + \frac{\Delta x^3 (-2 \I p_0 \omega ^3 (s-\sigma )+6
   p_0 \omega ^4 \epsilon +u_0 \omega ^3 (4 \omega 
   \epsilon -3 \I s)-3 \I \sigma  u_0 \omega ^3)}{288 s} \Big]^\top, \label{eq:jsdecomp2}
\end{align}
\begin{align*}
\tilde V_3 &= \left[ 0, ~0, ~\frac{(6\epsilon\omega\hat{p}_0 + \I(\tilde{s}-\sigma)\hat{u}_0)\omega^3\Delta x^3}{144\tilde{s}}, ~\frac{(6\epsilon\omega\hat{u}_0 + \I(\tilde{s}+\sigma)\hat{p}_0)\omega^3\Delta x^3}{144\tilde{s}} \right]^\top, \nonumber \\
\tilde V_4 &= \left[ 0, ~0, ~\frac{(-6\epsilon\omega\hat{p}_0 + \I(\tilde{s}+\sigma)\hat{u}_0)\omega^3\Delta x^3}{144\tilde{s}}, ~\frac{(-6\epsilon\omega\hat{u}_0 + \I(\tilde{s}-\sigma)\hat{p}_0)\omega^3\Delta x^3}{144\tilde{s}} \right]^\top.\nonumber
\end{align*}

In this decomposition, consider first the contribution of the spurious eigenvectors $\tilde V_3, \tilde V_4$. It is of order $\mathcal O(\Delta x^3)$, i.e., of the order of accuracy of the method, uniformly in $\epsilon$: expanding the stated terms in $\tilde V_3$ and $\tilde V_4$ as a power series in $\epsilon$, one finds
\begin{equation*}
 \left[ 0, ~0, ~\mathcal O(\epsilon), ~\frac{1}{72} \I \Delta x^3 \omega^3 \hat p_0 + \mathcal O(\epsilon)\right], \quad
 \left[ 0, ~0, ~\frac{1}{72} \I \Delta x^3 \omega^3 \hat u_0 + \mathcal O(\epsilon), ~\mathcal O(\epsilon)\right],
\end{equation*}
up to terms $\mathcal{O}(\Delta x^4)$.
The decomposition along the two physical modes $\tilde V_1,\tilde V_2$ quantifies which part of the data is evolving according to the first (physical) eigenvalue $\tilde \lambda_1 \simeq \lambda_1$, and which - according to the second. In order to determine the error, one needs to compare \eqref{eq:jsdecomp1}-\eqref{eq:jsdecomp2} to the averages and point values of the same decomposition at the PDE level. These latter are obtained by replacing $\hat p_0$ and $\hat u_0$ in \eqref{eq:projfouriermodeontoDof} by the respective components of $V_1$ and $V_2$ from \eqref{eq:eigenvectorpde1}--\eqref{eq:eigenvectorpde2}. After subtraction, the error of $\tilde V_1$ (the mode corresponding to the limit heat equation) reads
\begin{equation*}
\Big [0, ~0, 
~\frac{(2\epsilon\omega\hat{p}_0+\I(s-\sigma)\hat{u}_0)\omega^3\Delta x^3}{144 s},
 ~\frac{(2\epsilon\omega\hat{u}_0+\I(s+\sigma)\hat{p}_0)\omega^3\Delta x^3}{144 s} \Big],
\end{equation*}
where $\mathcal{O}(\Delta x^4)$ terms are omitted.
The third entry is a term $O(\epsilon)$, while the fourth is $\frac{1}{72} \I\Delta x^3 \omega^3 \hat p_0 + \mathcal O(\epsilon)$.
The error of $\tilde{V}_2$ is
\begin{equation*}
    \Big[
    0,~
    0,
    ~\frac{(-2\epsilon\omega\hat{p}_0+\I(s+\sigma)\hat{u}_0)\omega^3\Delta x^3}{144 s},~
    \frac{(-2\epsilon\omega\hat{u}_0+\I(s-\sigma)\hat{p}_0)\omega^3\Delta x^3}{144 s}\Big],
\end{equation*}
where the third entry is $\frac{1}{72} \I\Delta x^3 \omega^3 \hat u_0 + \mathcal O(\epsilon)$, and the fourth entry is $\mathcal{O}(\epsilon)$.
Thus, these errors are of third order uniformly in $\epsilon$.

The analysis thus shows that, in principle, the relevant errors of the Fourier decomposition are $\mathcal O(\Delta x^3)$. For the eigenvectors, this is uniformly in $\epsilon$. Before condensing these results into a statement on the error of the method, we perform analogous computations for the point value update based on the alternating flux. With the two sets of results in place, in \Cref{ssec:discussion} we present the remainder of the error analysis based on the Fourier decomposition.

\subsubsection{Point value update using the alternating flux}

The same computations are now performed for the AF method with the point update based on the prescription from \Cref{sec:alternating}. As mentioned there, this kind of upwinding does not correspond to the usual way of achieving stability. However, it is stable in one spatial dimension and the comparison between the two kinds of point value updates helps elucidating the behavior of the numerical methods in the limit $\epsilon \to 0$.

When the point values are updated with the ``alternating flux'' prescription, the scheme is
\begin{equation*}
	\frac{\dd}{\dd t} [\hat{p}_1, \hat{u}_1, \hat{p}_2, \hat{u}_2]^\top
	= \bm{G} [\hat{p}_1, \hat{u}_1, \hat{p}_2, \hat{u}_2]^\top, 
\end{equation*}
with the evolution matrix
\begin{equation*}
	\bm{G}
	= \frac{1}{\epsilon\Delta x}\begin{bmatrix}
		0 & 0 & 0 & t_x^{-1}-1 \\
		0 & -\epsilon^{-1}\sigma\Delta x & t_x^{-1}-1 & 0 \\
		0 & 6 & 0 & -4-2t_x^{-1} \\
		 -6t_x & 0 & 4+2t_x & -\epsilon^{-1}\sigma\Delta x \\
	\end{bmatrix}.
\end{equation*}

It can be diagonalized explicitly.
The eigenvalues are 
\begin{align*}
    \tilde{\lambda}_1 &= \lambda_1  + \frac{\omega^6 \Delta x^4}{540 s}  + \mathcal{O}(\Delta x^5), \\
    \tilde{\lambda}_2 & = \lambda_2  - \frac{\omega^6 \Delta x^4}{540 s}  + \mathcal{O}(\Delta x^5), \\
    \tilde{\lambda}_3 &= -\frac{6\I}{\epsilon\Delta x} - \frac{\sigma}{2\epsilon^2} + \I \Delta x \frac{12 \epsilon^2 \omega^2 + \sigma^2}{48 \epsilon^3} + \mathcal{O}(\Delta x^3), \\
    \tilde{\lambda}_4 &= \frac{6\I}{\epsilon\Delta x} - \frac{\sigma}{2\epsilon^2} - \I \Delta x \frac{12 \epsilon^2 \omega^2 + \sigma^2}{48 \epsilon^3} + \mathcal{O}(\Delta x^3).
\end{align*}
Observe first of all that this time, the physical eigenvalues $\tilde \lambda_1, \tilde \lambda_2$ are 4th-order accurate, and that the error is uniform in $\epsilon$. This is a difference from the results with the JS obtained above. We refrain from stating the entire expansion of the spurious eigenvalues $\tilde \lambda_3, \tilde \lambda_4$, since the contribution of the spurious modes is again negligible: Upon decomposing the initial data $\hat w_0$ in the eigenvectors, the terms
$\tilde V_3$ and $\tilde V_4$ are, neglecting contributions $\mathcal O(\Delta x^4)$,
\begin{align*}
\tilde V_3 &= \left[ 0, ~0, ~\frac{1}{144} \Delta x^3 \omega^3 (\hat u_0 - \I \hat p_0), ~\frac{1}{144} \Delta x^3 \omega^3 (\hat p_0 + \I \hat u_0)\right ], \\
\tilde V_4 &= \left[ 0, ~0, -\frac{1}{144} \I \Delta x^3 \omega^3 (\hat p_0 - \I \hat u_0), -\frac{1}{144} \Delta x^3 \omega^3 (\hat p_0 - \I \hat u_0)\right ].
\end{align*}
Observe that these are independent of $\epsilon$.
Finally, we compute the difference between the projection of the exact mode corresponding to $\lambda_1$ onto the degrees of freedom and the corresponding discrete mode $\tilde V_1$, finding
\begin{equation*}
 \left[ \frac{\epsilon \omega^4 \hat u_0 \Delta x^3}{72 s}, -\frac{\epsilon \omega^4 \hat p_0 \Delta x^3}{72 s}, -\frac{1}{144} \I \omega^3 \left(1 + \frac{\sigma}{s} \right ) \hat p_0\Delta x^3, \frac{1}{144} \I \omega^3 \left(1 - \frac{\sigma}{s} \right ) \hat u_0 \Delta x^3\right ]
\end{equation*}
up to terms $\mathcal O(\Delta x^4$).
The error of the mode corresponding to the initial layer is
\begin{equation*}
 \left[ -\frac{\epsilon \omega^4 \hat u_0 \Delta x^3}{72 s}, \frac{\epsilon \omega^4 \hat p_0 \Delta x^3}{72 s}, -\frac{1}{144} \I \omega^3 \left(1 - \frac{\sigma}{s} \right ) \hat p_0\Delta x^3, \frac{1}{144} \I \omega^3 \left(1 + \frac{\sigma}{s} \right ) \hat u_0 \Delta x^3\right ]
\end{equation*}
up to terms $\mathcal O(\Delta x^4$).
All terms are of the order of accuracy of the method and also $\mathcal O(\epsilon)$.

\subsubsection{Discussion} \label{ssec:discussion}

Following \cite{Guo_2013_JoCP_SuperconvergenceDiscontinuous} we observe that the error of the method can be related to the errors of the eigendecomposition of $\bm{G}$ as performed in the previous Sections in the following way: 
\begin{align*}
 \left \| \pi \hat q(t) - \hat q_h(t) \right \|_{L^2} &= \left \| \sum_{k=1}^2 (\pi V_k \exp(\lambda_k t) - \tilde V_k \exp(\tilde \lambda_k t))  + \sum_{k=3}^4 \tilde V_k \exp(\tilde \lambda_k t)) \right \|_{L^2} \\
 &\leq  \sum_{k=1}^2 \| \pi V_k - \tilde V_k \|  | \exp(\lambda_k t) | + \sum_{k=1}^2  \| \tilde V_k \| | \exp(\lambda_k t) - \exp(\tilde \lambda_k t)|  \\ & \qquad  \nonumber + \sum_{k=3}^4 \| \tilde V_k \| | \exp(\tilde \lambda_k t) |  .
\end{align*}
Here $\pi$ denotes the projection \eqref{eq:projfouriermodeontoDof} of a mode $[\hat p_0, \hat u_0]$ associated to the PDE (a vector in $\mathbb C^2$) onto a Fourier mode associated to the degrees of freedom of the scheme (a vector in $\mathbb C^4$).

For both the JS and the alternating flux, we have seen that the contribution of the spurious modes is negligible, i.e., $\| \tilde V_k \| \in \mathcal O(\Delta x^3)$ for $k=3,4$. Similarly, the error $\| \pi V_k - \tilde V_k \|$, $k=1,2$ of the physical modes is $\mathcal O(\Delta x^3)$ as well.

Finally, 
\begin{equation*}
| \exp(\lambda_k t) - \exp(\tilde \lambda_k t)| \leq C t |\lambda_k - \tilde \lambda_k|.
\end{equation*}
This value is $\mathcal O(\Delta x^4)$ for the alternating flux, and $\mathcal O\left(\frac{\Delta x^3}{\epsilon} \right)$ for the JS. We conclude that in the limit, the scheme based on JS suffers a loss of accuracy, while the one based on the ``alternating flux'' prescription does not.
The loss of accuracy can be alleviated by choosing $\Delta x$ sufficiently small in relation to $\epsilon$. This is confirmed in Example \ref{ex:1d_accuracy}.
The formal analysis presented here, however, is insufficient to obtain theoretical convergence rates, which is beyond the scope of the paper. Below, we present detailed numerical studies of the order of accuracy for various values of $\epsilon$.

\section{2D active flux scheme}\label{sec:2d_af_scheme}
This section presents the 2D semi-discrete AF methods for the hyperbolic heat equation
\begin{equation*}
	\left\{~
	\begin{aligned}
		&p_t + \frac{1}{\epsilon}u_x + \frac{1}{\epsilon}v_y = 0, \\
		&u_t + \frac{1}{\epsilon}p_x + \frac{\sigma}{\epsilon^2}u = 0, \\
		&v_t + \frac{1}{\epsilon}p_y + \frac{\sigma}{\epsilon^2}v = 0.
	\end{aligned}
	\right.
\end{equation*}

Consider a 2D uniform Cartesian mesh with $N_1\times N_2$ cells
$I_{ij} = [x_{\xl}, x_{\xr}]\times [y_{\yl}, y_{\yr}]$ with cell centers $((x_{\xl} + x_{\xr})/2, (y_{\yl} + y_{\yr})/2)$ and cell sizes $\Delta x, \Delta y$.
The degrees of freedom contain the cell averages, face-centered point values, and nodal values at corners.
The cell average is updated following the finite volume manner
\begin{equation}\label{eq:2d_av}
	\begin{aligned}
		\frac{\dd \bar{p}_{ij}}{\dd t} =& - \frac{1}{\epsilon\Delta x}\left(\hat{u}_{\xr,j} - \hat{u}_{\xl,j}\right)
		- \frac{1}{\epsilon\Delta y}\left(\hat{v}_{i,\yr} - \hat{v}_{i,\yl}\right), \\
		\frac{\dd \bar{u}_{ij}}{\dd t} =& - \frac{1}{\epsilon\Delta x}\left(\hat{p}_{\xr,j} - \hat{p}_{\xl,j}\right) - \frac{\sigma}{\epsilon^2}\bar{u}_{ij}, \\
		\frac{\dd \bar{v}_{ij}}{\dd t} =& - \frac{1}{\epsilon\Delta y}\left(\hat{p}_{i,\yr} - \hat{p}_{i,\yl}\right) - \frac{\sigma}{\epsilon^2}\bar{v}_{ij},
	\end{aligned}
\end{equation}
where the numerical fluxes are
\begin{align*}
	&\hat{q}_{\xr,j} = \frac{1}{6}\left(q_{\xr,\yl} + 4q_{\xr,j} + q_{\xr,\yr}\right), \\
	&\hat{q}_{i,\yr} = \frac{1}{6}\left(q_{\xl,\yr} + 4q_{i,\yr} + q_{\xr,\yr}\right), \quad q = p, u, v.
\end{align*}

For the point value update, the Jacobian in the $x$- and $y$-directions are split as
\begin{equation*}
	\bm{J}_x^+ = \frac{1}{\epsilon}\begin{pmatrix}
		\frac12 & \frac12 & 0 \\ \frac12 & \frac12 & 0 \\ 0 & 0 & 0 \\
	\end{pmatrix},\quad
	\bm{J}_x^- = \frac{1}{\epsilon}\begin{pmatrix}
		-\frac12 & \frac12 & 0 \\ \frac12 & -\frac12 & 0 \\ 0 & 0 & 0 \\
	\end{pmatrix},
\end{equation*}
\begin{equation*}
	\bm{J}_y^+ = \frac{1}{\epsilon}\begin{pmatrix}
		\frac12 & 0 & \frac12 \\ 0 & 0 & 0 \\ \frac12 & 0 & \frac12 \\
	\end{pmatrix},\quad
	\bm{J}_y^- = \frac{1}{\epsilon}\begin{pmatrix}
		-\frac12 & 0 & \frac12 \\ 0 & 0 & 0 \\ \frac12 & 0 & -\frac12 \\
	\end{pmatrix}.
\end{equation*}
Then using the JS-based discretizations in \cite{Duan_2025_Active_SJoSC,Abgrall_2025_Semi_JoSC}, the point value update at the corner is
\begin{equation}\label{eq:2d_pnt_node}
	\begin{aligned}
		\frac{\dd p_{\xr,\yr}}{\dd t} =& - \frac{1}{2\epsilon}\Big[ (D_x^+)_{\xr,\yr}(p+u) - (D_x^-)_{\xr,\yr}(p-u)\Big] \\
		&- \frac{1}{2\epsilon}\Big[ (D_y^+)_{\xr,\yr}(p+v) - (D_y^-)_{\xr,\yr}(p-v)\Big], \\
		\frac{\dd u_{\xr,\yr}}{\dd t} =& - \frac{1}{2\epsilon}\Big[(D_x^+)_{\xr,\yr}(p+u) + (D_x^-)_{\xr,\yr}(p-u) \Big] - \frac{\sigma}{\epsilon^2} u_{i+\frac12,\yr}, \\
		\frac{\dd v_{\xr,\yr}}{\dd t} =& - \frac{1}{2\epsilon}\Big[(D_y^+)_{\xr,\yr}(p+v) + (D_y^-)_{\xr,\yr}(p-v) \Big] - \frac{\sigma}{\epsilon^2} v_{i+\frac12,\yr},
	\end{aligned}
\end{equation}
where
\begin{align*}
	&(D_x^+)_{\xr,\yr}(a) = a_{\xl,\yr} - 4a_{i,\yr} + 3a_{\xr,\yr}, \\
	&(D_x^-)_{\xr,\yr}(a) = -3a_{\xr,\yr} + 4a_{i+1,\yr} - a_{i+\frac32,\yr}.
\end{align*}
The point value update at the face center $(x_{\xr},y_j)$ is
\begin{equation}\label{eq:2d_pnt_facex}
	\begin{aligned}
		\frac{\dd p_{\xr,j}}{\dd t} =& - \frac{1}{2\epsilon}\Big[ (D_x^+)_{\xr,j}(p+u) - (D_x^-)_{\xr,j}(p-u)\Big] - \frac{1}{\epsilon} (D_y)_{\xr,j}v, \\
		\frac{\dd u_{\xr,j}}{\dd t} =& - \frac{1}{2\epsilon}\Big[(D_x^+)_{\xr,j}(p+u) + (D_x^-)_{\xr,j}(p-u) \Big] - \frac{\sigma}{\epsilon^2} u_{i+\frac12,j}, \\
		\frac{\dd v_{\xr,j}}{\dd t} =& - \frac{1}{\epsilon}(D_y)_{\xr,j}p - \frac{\sigma}{\epsilon^2} v_{i+\frac12,j},
	\end{aligned}
\end{equation}
where
\begin{align*}
	&(D_x^+)_{\xr,j}(a) = a_{\xl,j} - 4a_{i,j} + 3a_{\xr,j}, \\
	&(D_x^-)_{\xr,j}(a) = -3a_{\xr,j} + 4a_{i+1,j} - a_{i+\frac32,j}, \\
	&(D_y)_{\xr,j}(a) = a_{\xr,\yr} - a_{\xr,\yl},
\end{align*}
and the cell center value is obtained by
\begin{align}
	a_{i,j} =&\ \frac{1}{16}\Big[36\bar{a}_{i,j} - 4\left(a_{\xl,j} + a_{\xr,j} + a_{i,\yl} + a_{i,\yr}\right) \nonumber\\
	&\qquad - \left(a_{\xl,\yl} + a_{\xr,\yl} + a_{\xl,\yr} + a_{\xr,\yr}\right) \Big]. \label{eq:2d_cell_center}
\end{align}
The point value update at the face center $(x_i, y_{\yr})$ is
\begin{equation}\label{eq:2d_pnt_facey}
	\begin{aligned}
		\frac{\dd p_{i,\yr}}{\dd t} =& - \frac{1}{\epsilon} (D_x)_{i,\yr}u - \frac{1}{2\epsilon}\Big[ (D_y^+)_{i,\yr}(p+v) - (D_y^-)_{i,\yr}(p-v)\Big], \\
		\frac{\dd u_{i,\yr}}{\dd t} =& - \frac{1}{\epsilon}(D_x)_{i,\yr}p - \frac{\sigma}{\epsilon^2} u_{i,\yr}, \\
		\frac{\dd v_{i,\yr}}{\dd t} =& - \frac{1}{2\epsilon}\Big[(D_y^+)_{i,\yr}(p+v) + (D_y^-)_{i,\yr}(p-v) \Big] - \frac{\sigma}{\epsilon^2} v_{i,\yr}.
	\end{aligned}
\end{equation}

\subsection{Formal asymptotic analysis}
Similar to the analysis in Section \ref{sec:1d_af_formal}, under the assumption of sufficient regularity of numerical solutions, one can first match the terms of order $\epsilon^{-2}$ in \eqref{eq:2d_av}-\eqref{eq:2d_pnt_facey} and obtain
\begin{equation}\label{eq:2d_limit_u0}
	\bar{u}_{ij}^{(0)} = \bar{v}_{ij}^{(0)} = u_{\zeta}^{(0)} = v_{\zeta}^{(0)} = 0,
\end{equation}
where $\zeta$ denotes the point values at face centers or corners.

For the terms of order $\epsilon^{-1}$, one has
\begin{equation}\label{eq:2d_limit_u1}
	\begin{aligned}
		&\bar{u}_{ij}^{(1)} = - \frac{1}{\sigma\Delta x}\left(\hat{p}_{\xr,j}^{(0)} - \hat{p}_{\xl,j}^{(0)}\right),~
		\bar{v}_{ij}^{(1)} = - \frac{1}{\sigma\Delta y}\left(\hat{p}_{i,\yr}^{(0)} - \hat{p}_{i,\yl}^{(0)}\right), \\
		&u_{\xr,\yr}^{(1)} = - \frac{1}{2\sigma}\Big[(D_x^+)_{\xr,\yr} + (D_x^-)_{\xr,\yr} \Big]p^{(0)}, \\
		&v_{\xr,\yr}^{(1)} = - \frac{1}{2\sigma}\Big[(D_y^+)_{\xr,\yr} + (D_y^-)_{\xr,\yr}\Big]p^{(0)}, \\
		&u_{\xr,j}^{(1)} = - \frac{1}{2\sigma}\Big[(D_x^+)_{\xr,j} + (D_x^-)_{\xr,j} \Big]p^{(0)}, \\
		&v_{\xr,j}^{(1)} = - \frac{1}{\sigma}(D_y)_{\xr,j}p^{(0)}, \\
		&u_{i,\yr}^{(1)} = - \frac{1}{\sigma}(D_x)_{i,\yr}p^{(0)}, \\
		&v_{i,\yr}^{(1)} = - \frac{1}{2\sigma}\Big[(D_y^+)_{i,\yr} + (D_y^-)_{i,\yr} \Big]p^{(0)},
	\end{aligned}
\end{equation}
which are generally some approximations of $u^{(1)} = -\frac{1}{\sigma}p^{(0)}_x, v^{(1)} = -\frac{1}{\sigma}p^{(0)}_y$.
Take the terms related to $u^{(1)}$ as examples while those of $v^{(1)}$ can be treated similarly.
The Simpson's rule is 4th order,
\begin{equation*}
	\hat{p}^{(0)}_{\xr,j} = \frac{1}{6}\left( p^{(0)}_{\xr,\yl} + 4p^{(0)}_{\xr,j} + p^{(0)}_{\xr,\yr} \right) 
	= \left(\frac{1}{\Delta y}\int_{y_{\yl}}^{y_{\yr}} p^{(0)} \dd y\right)\Bigg|_{x_{\xr}} + \mathcal{O}(\Delta y^4),
\end{equation*}
thus
\begin{align*}
	\bar{u}_{ij}^{(1)} =& - \frac{1}{\sigma\Delta x}\left(\frac{1}{\Delta y}\int_{y_{\yl}}^{y_{\yr}} p^{(0)} \dd y\right)\Bigg|_{x_{\xl}}^{x_{\xr}}
	+ \mathcal{O}(\Delta y^4),
\end{align*}
which is a finite volume approximation of $u^{(1)} = -\frac{1}{\sigma}p^{(0)}_x$.
By using the Taylor expansion,
\begin{align*}
	u_{i+\frac12,\yr}^{(1)} =& - \frac{1}{2\sigma}\Big[(D_x^+)_{\xr,\yr} + (D_x^-)_{\xr,\yr} \Big]p^{(0)}
	= -\frac{1}{\sigma}p^{(0)}_x \Big|_{x_{\xr},y_{\yr}} + \mathcal{O}(\Delta x^2), \\
	u_{\xr,j}^{(1)} =& - \frac{1}{2\sigma}\Big[(D_x^+)_{\xr,j} + (D_x^-)_{\xr,j} \Big]p^{(0)}
	= -\frac{1}{\sigma}p^{(0)}_x \Big|_{x_{\xr},y_{j}} + \mathcal{O}(\Delta x^2), \\
	u_{i,\yr}^{(1)} =& - \frac{1}{\sigma}(D_x)_{i,\yr}p^{(0)}
	= -\frac{1}{\sigma}p^{(0)}_x \Big|_{x_{i},y_{\yr}} + \mathcal{O}(\Delta x^2),
\end{align*}
so that they are approximations of $u^{(1)} = -\frac{1}{\sigma}p^{(0)}_x$ at corresponding points.

Collecting the terms of order $\epsilon^0$ in \eqref{eq:2d_av}-\eqref{eq:2d_pnt_facey} gives
\begin{align}
	\frac{\dd}{\dd t}\bar{p}_{ij}^{(0)} =&
	- \frac{1}{\Delta x}\left(\hat{u}_{\xr,j}^{(1)} - \hat{u}_{\xl,j}^{(1)}\right)
	- \frac{1}{\Delta y}\left(\hat{v}_{i,\yr}^{(1)} - \hat{v}_{i,\yl}^{(1)}\right), \label{eq:2d_limit_av_p0} \\
	\frac{\dd p_{\xr,\yr}^{(0)}}{\dd t} =& - \frac{1}{2}\Big[ (D_x^+)_{\xr,\yr}(p^{(1)}+u^{(1)}) - (D_x^-)_{\xr,\yr}(p^{(1)}-u^{(1)})\Big] \nonumber\\
	& - \frac{1}{2}\Big[ (D_y^+)_{\xr,\yr}(p^{(1)}+v^{(1)}) - (D_y^-)_{\xr,\yr}(p^{(1)}-v^{(1)})\Big], \nonumber\\
	\frac{\dd p_{\xr,j}^{(0)}}{\dd t} =& - \frac{1}{2}\Big[ (D_x^+)_{\xr,j}(p^{(1)}+u^{(1)}) - (D_x^-)_{\xr,j}(p^{(1)}-u^{(1)})\Big] - (D_y)_{\xr,j}(v^{(1)}), \nonumber\\
	\frac{\dd p_{i,\yr}^{(0)}}{\dd t} =& - (D_x)_{i,\yr}u^{(1)} - \frac{1}{2}\Big[ (D_y^+)_{i,\yr}(p^{(1)}+v^{(1)}) - (D_y^-)_{i,\yr}(p^{(1)}-v^{(1)})\Big], \nonumber
\end{align}
where the last three can be simplified based on the smoothness of $p^{(1)}$ and the Taylor expansion as
\begin{equation}\label{eq:2d_limit_p0}
	\begin{aligned}
		\frac{\dd p_{\xr,\yr}^{(0)}}{\dd t} =& - \frac{1}{2}\Big[ (D_x^+)_{\xr,\yr} + (D_x^-)_{\xr,\yr}\Big]u^{(1)} 
		- \frac{1}{2}\Big[ (D_y^+)_{\xr,\yr} + (D_y^-)_{\xr,\yr}\Big]v^{(1)} \\
		&+\mathcal{O}(\Delta x^3 + \Delta y^3), \\
		\frac{\dd p_{\xr,j}^{(0)}}{\dd t} =& - \frac{1}{2}\Big[ (D_x^+)_{\xr,j} + (D_x^-)_{\xr,j}\Big]u^{(1)} - (D_y)_{\xr,j}v^{(1)}
		+\mathcal{O}(\Delta x^3), \\
		\frac{\dd p_{i,\yr}^{(0)}}{\dd t} =& - (D_x)_{i,\yr}u^{(1)} - \frac{1}{2}\Big[ (D_y^+)_{i,\yr} + (D_y^-)_{i,\yr}\Big]v^{(1)}
		+\mathcal{O}(\Delta y^3).
	\end{aligned}
\end{equation}
Note that $\bar{u}_{ij}^{(1)}, u_{\zeta}^{(1)}, \bar{v}_{ij}^{(1)}, v_{\zeta}^{(1)}$ appeared in \eqref{eq:2d_limit_p0} are given in \eqref{eq:2d_limit_u1},
and the cell center value is given by \eqref{eq:2d_cell_center}.

One can further examine the accuracy of \eqref{eq:2d_limit_av_p0}.
Since
\begin{align*}
	- \hat{u}^{(1)}_{\xr,j} =& -\frac{1}{6}\left( u^{(1)}_{\xr,\yl} + 4u^{(1)}_{\xr,j} + u^{(1)}_{\xr,\yr} \right) \\
	=&\ \frac{1}{12\sigma}\Big[(D_x^+)_{\xr,\yl} + (D_x^-)_{\xr,\yl}
	+ 4(D_x^+)_{\xr,j} \\
	&+ 4(D_x^-)_{\xr,j}
	+ (D_x^+)_{\xr,\yr} + (D_x^-)_{\xr,\yr} \Big]p^{(0)} \\
	=&\ \frac{1}{\sigma}\left(\frac{1}{\Delta y}\int_{y_{\yl}}^{y_{\yr}} p^{(0)}_x \dd y\right)\Bigg|_{x_{\xr}} + \mathcal{O}(\Delta x^2 + \Delta y^4),
\end{align*}
treating other numerical fluxes similarly, then \eqref{eq:2d_limit_av_p0} becomes
\begin{equation}\label{eq:2d_av_limit_appro}
	\frac{\dd}{\dd t}\bar{p}_{ij}^{(0)} = \frac{1}{\sigma\Delta x}\left(\frac{1}{\Delta y}\int_{y_{\yl}}^{y_{\yr}} p^{(0)}_x \dd y\right)\Bigg|_{x_{\xl}}^{x_{\xr}}
	+ \frac{1}{\sigma\Delta y}\left(\frac{1}{\Delta x}\int_{x_{\xl}}^{x_{\xl}} p^{(0)}_y \dd x\right)\Bigg|_{y_{\yl}}^{y_{\yr}}
	+ \mathcal{O}(\Delta x^2 + \Delta y^2).
\end{equation}
Thus it is a finite volume approximation of $p^{(0)}_t = \frac{1}{\sigma}p^{(0)}_{xx} + \frac{1}{\sigma}p^{(0)}_{yy}$ with the truncation error $\mathcal{O}(\Delta x^2 + \Delta y^2)$.
By the Taylor expansion, the other equations in \eqref{eq:2d_limit_p0} can be written as
\begin{equation}\label{eq:2d_pnt_limit_appro}
		\frac{\dd p_{\zeta}^{(0)}}{\dd t} = \left(\frac{1}{\sigma}p^{(0)}_{xx} + \frac{1}{\sigma}p^{(0)}_{yy}\right)\Bigg|_{\zeta} + \mathcal{O}(\Delta x^2 + \Delta y^2),
\end{equation}
for all the point values at face centers or corners.

\begin{proposition}
	Assuming sufficient regularity of the solutions, as $\epsilon\rightarrow 0$, the leading order solutions of the 2D AF schemes \eqref{eq:2d_av}-\eqref{eq:2d_pnt_facey} satisfy \eqref{eq:2d_limit_u0}-\eqref{eq:2d_limit_p0},
	which are approximations of $p^{(0)}_t = \frac{1}{\sigma}p^{(0)}_{xx} + \frac{1}{\sigma}p^{(0)}_{yy}$ with truncation errors $\mathcal{O}(\Delta x^2 + \Delta y^2)$ in the sense of \eqref{eq:2d_av_limit_appro}-\eqref{eq:2d_pnt_limit_appro}.
\end{proposition}

\section{Numerical results}\label{sec:results}
This section conducts some numerical experiments to verify the convergence rates and AP property of the AF schemes.
Unless otherwise stated, periodic boundary conditions and the 3rd-order DIRK implemented in \cite{petsc-web-page} are used and the time step size is chosen as
\begin{equation*}
	\Delta t = C_{\texttt{CFL}} \min\{\Delta x, \Delta y\},
\end{equation*}
with the CFL number $C_\texttt{CFL} = 1.0$.
The linear system is solved by using GMRES with an incomplete LU factorization preconditioner \cite{petsc-web-page}.

\begin{example}[1D accuracy test]\label{ex:1d_accuracy}
	Consider the following exact solution \cite{Jang_2015_High_JoCP},
	\begin{equation*}
		p = \frac{1}{r}\exp(rt)\sin(x),~
		u = \epsilon\exp(rt)\cos(x),~
		r = \frac{-2}{1 + \sqrt{1-4\epsilon^2}},
	\end{equation*}
	on the domain $[0,2\pi]$ with periodic boundary conditions and $\sigma=1$.
	The test is solved until $T=1$ and the time step size is $0.2\Delta x^{4/3}$ to make the spatial error dominant.

	The errors and convergence rates are plotted in Figure \ref{fig:1d_accuracy} for $\epsilon=0.5,10^{-2},10^{-6}$,
	which shows that the JS-based AF scheme is 3rd-order accurate for $\epsilon=0.5$,
	while as $\epsilon\rightarrow 0$ the convergence rates reduce to $2$ except for the 4th-order convergence of the point value in $u$.
    When $\epsilon=10^{-2}$, it is verified that the scheme recovers the 3rd-order accuracy when the mesh size decreases to the magnitude of $\epsilon$. 
    The errors and convergence rates with the alternating flux are shown in Figure \ref{fig:1d_accuracy_alternating}.
    For $\epsilon=0.5$, one observes 3rd-order accuracy, while for $\epsilon=10^{-2}$ and $\epsilon=10^{-6}$, the convergence order is $4$ for the cell average of $p$ and point value of $u$, and $3$ for the other degrees of freedom.
    The appearance of the super-convergence (4th order instead of 3rd order) needs further study.

	\begin{figure}[htbp]
		\centering
		\begin{subfigure}[b]{0.325\textwidth}
			\centering
			\includegraphics[height=0.7\linewidth]{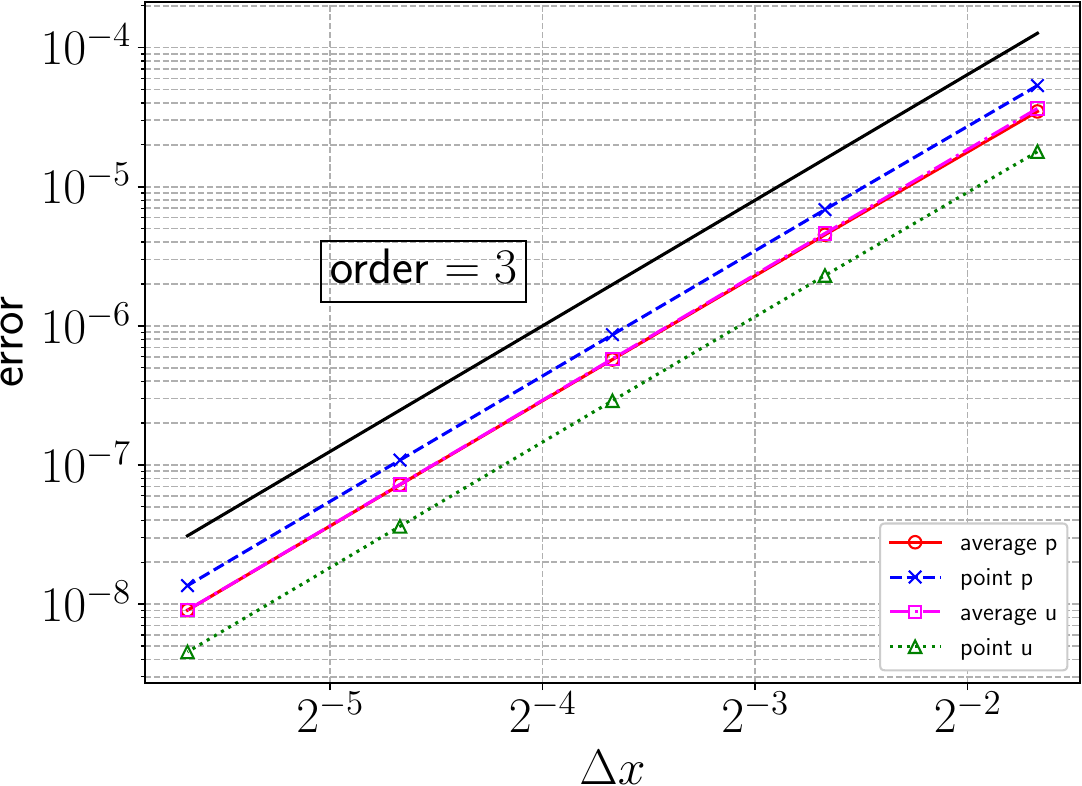}
		\end{subfigure}
		\begin{subfigure}[b]{0.325\textwidth}
			\centering
			\includegraphics[height=0.7\linewidth]{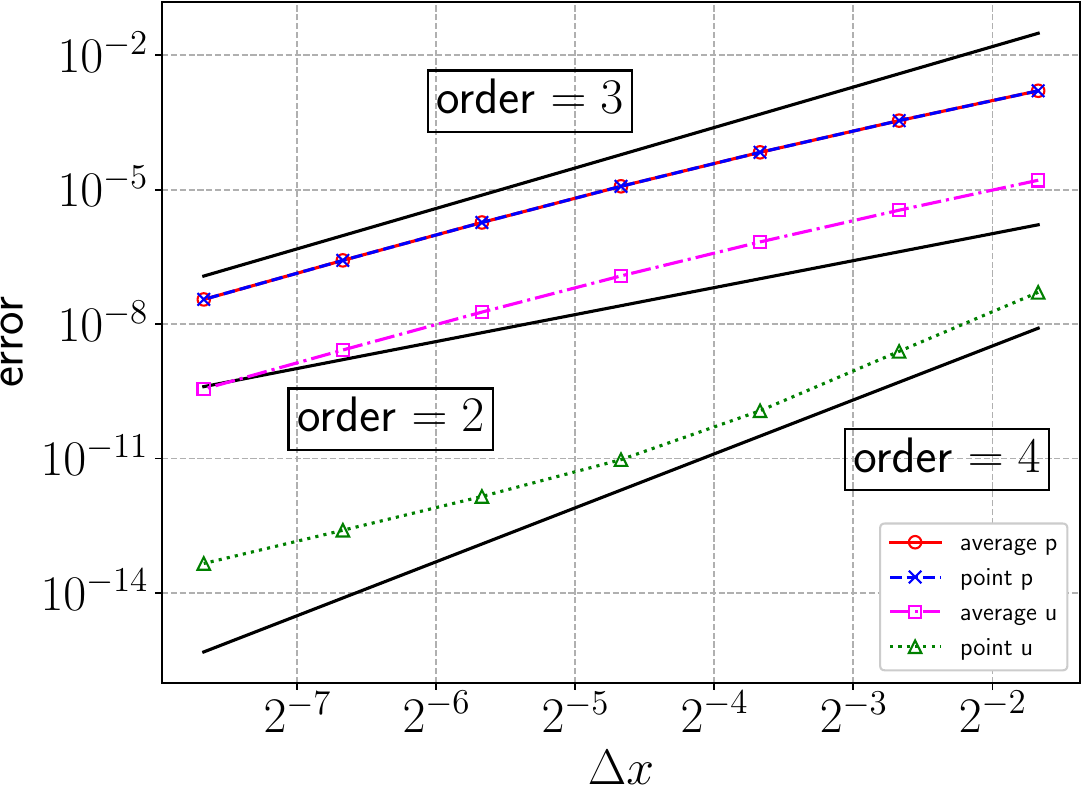}
		\end{subfigure}
		\begin{subfigure}[b]{0.325\textwidth}
			\centering
			\includegraphics[height=0.7\linewidth]{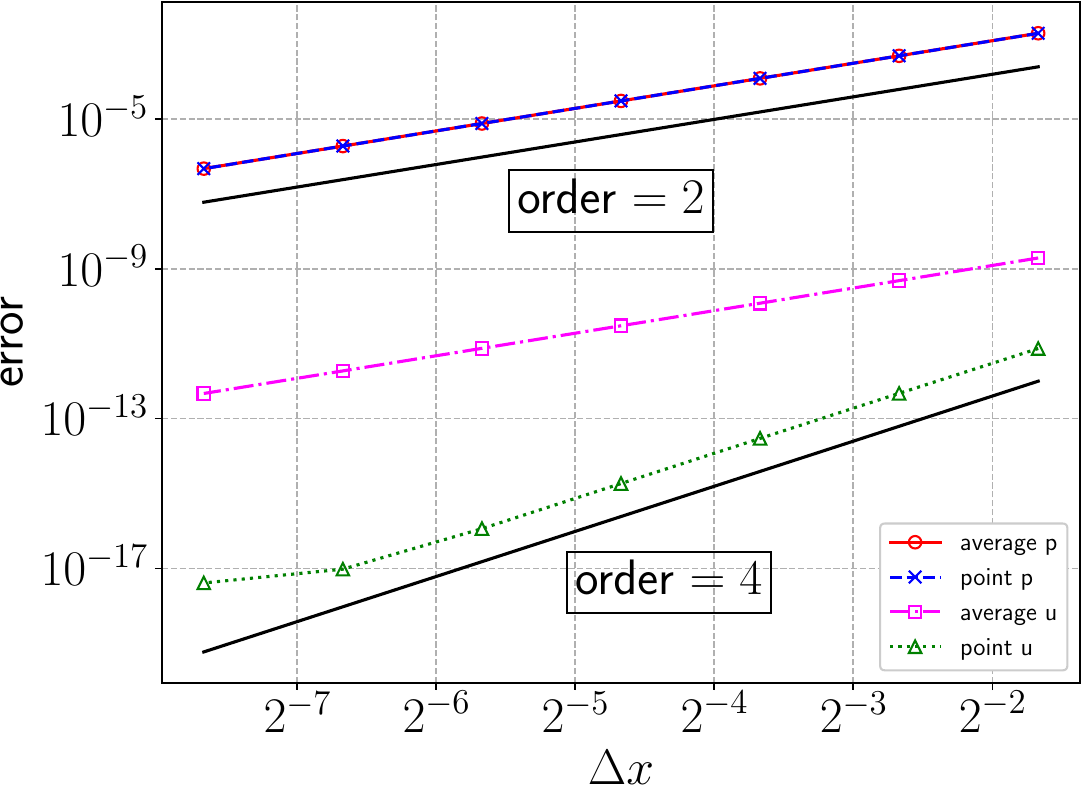}
		\end{subfigure}
		\caption{Example \ref{ex:1d_accuracy}.
			The errors and convergence rates obtained by the JS-based AF scheme and $\epsilon=0.5,10^{-2},10^{-6}$ from left to right.}
		\label{fig:1d_accuracy}
	\end{figure}

\begin{figure}[htbp]
		\centering
		\begin{subfigure}[b]{0.325\textwidth}
			\centering
			\includegraphics[height=0.7\linewidth]{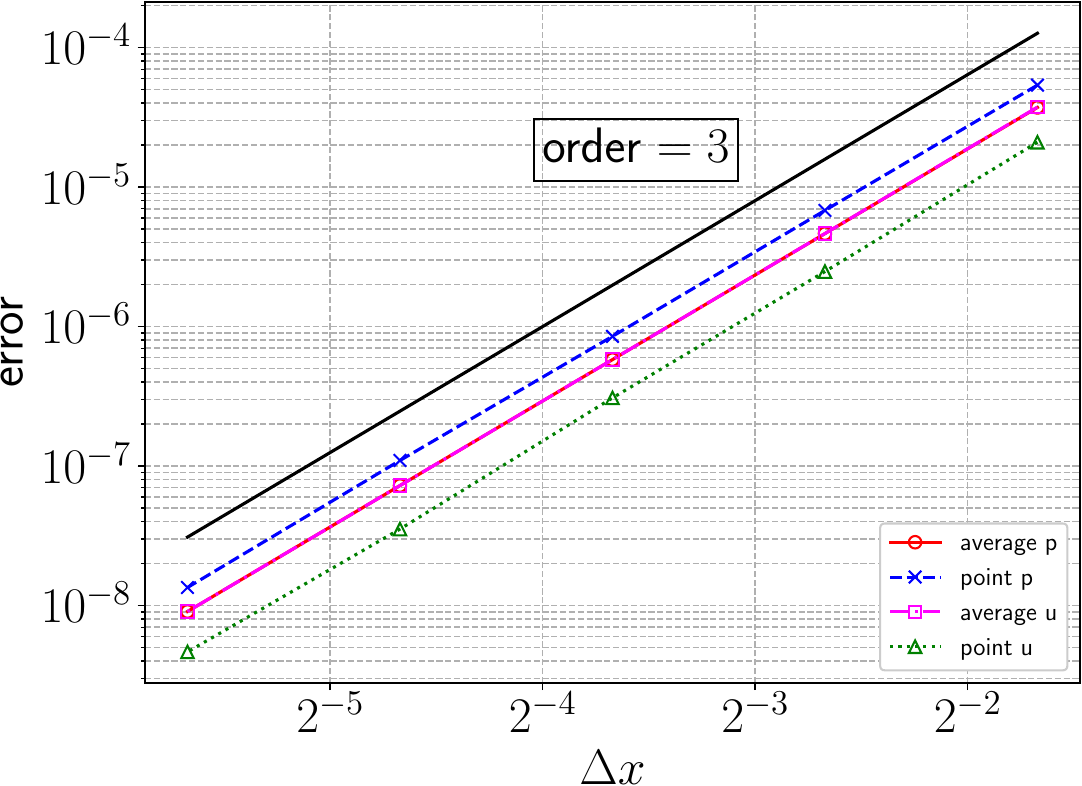}
		\end{subfigure}
		\begin{subfigure}[b]{0.325\textwidth}
			\centering
			\includegraphics[height=0.7\linewidth]{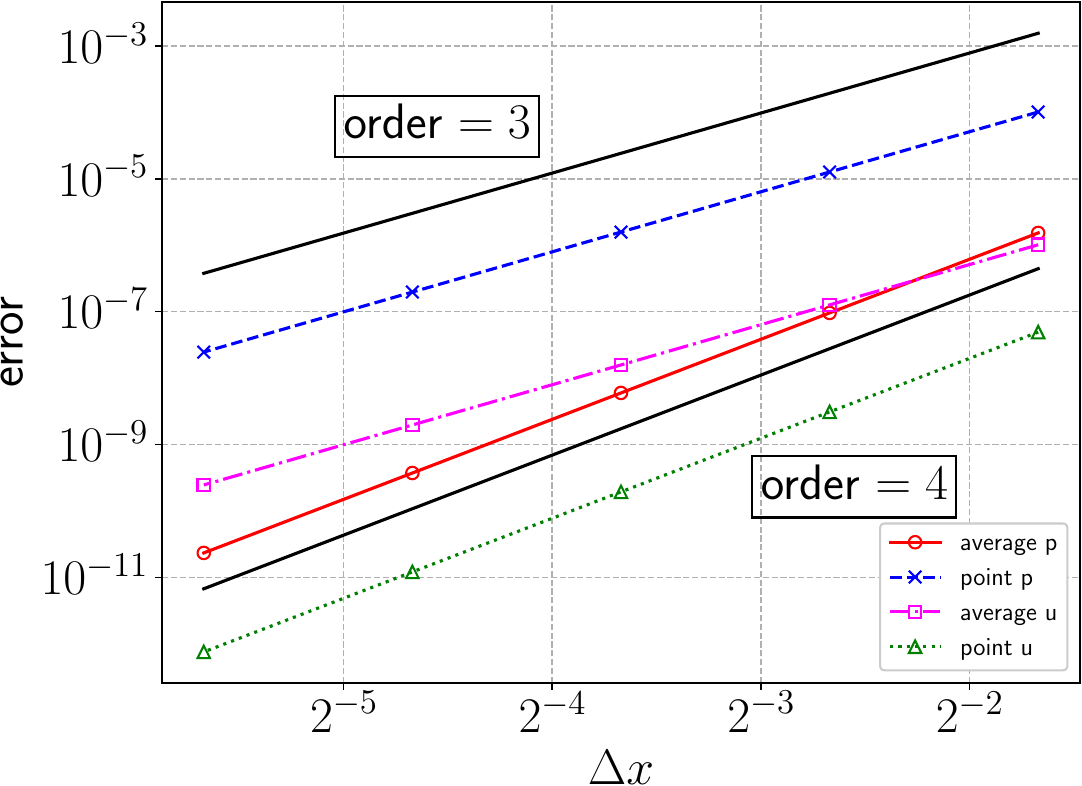}
		\end{subfigure}
		\begin{subfigure}[b]{0.325\textwidth}
			\centering
			\includegraphics[height=0.7\linewidth]{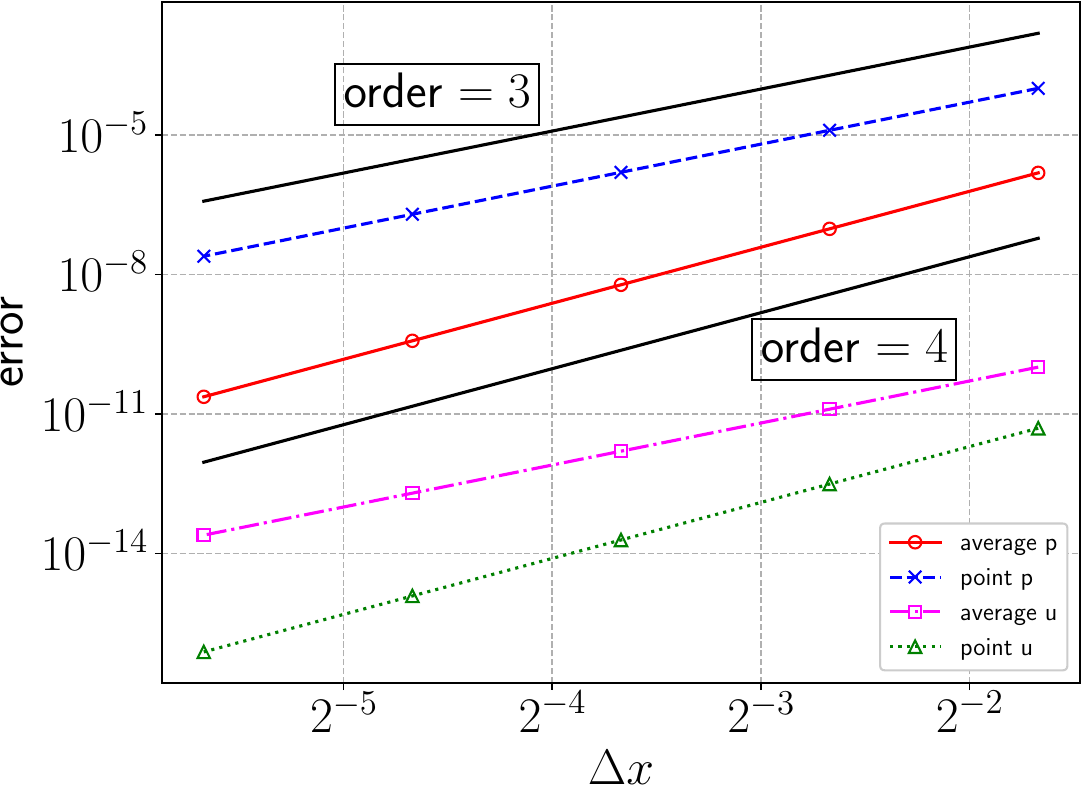}
		\end{subfigure}
		\caption{Example \ref{ex:1d_accuracy}.
			The errors and convergence rates obtained with the alternating flux and $\epsilon=0.5,10^{-2},10^{-6}$ from left to right.}
		\label{fig:1d_accuracy_alternating}
	\end{figure}
\end{example}

\begin{example}[Square wave]\label{ex:1d_square_wave}
	The initial $p$ of this test is a square wave, adapted from \cite{Jang_2015_High_JoCP}, 
	\begin{equation*}
		p = \begin{cases}
			2, &\text{if}~ \abs{x}<0.5, \\
			1, &\text{otherwise}, \\
		\end{cases}
	\end{equation*}
	and $u=0$ in the whole domain $[-1,1]$.
	Two cases are considered: the transport regime $\epsilon=0.7$ and the diffusive regime $\epsilon=10^{-6}$, with $\sigma=1$,
	and the final time is $0.25$ and $0.04$, respectively.
	
	The results obtained by the JS-based AF scheme on a coarse mesh of $40$ cells are shown in Figure \ref{fig:1d_square_wave}.
	It is observed that the scheme can capture shock waves for $\epsilon=0.7$ correctly,
	and can also capture the diffusive solution for $\epsilon=10^{-6}$ well.
	
	 \begin{figure}[htbp]
	 	\centering
	 	\begin{subfigure}[b]{0.45\textwidth}
	 		\centering
	 		\includegraphics[width=1.0\linewidth]{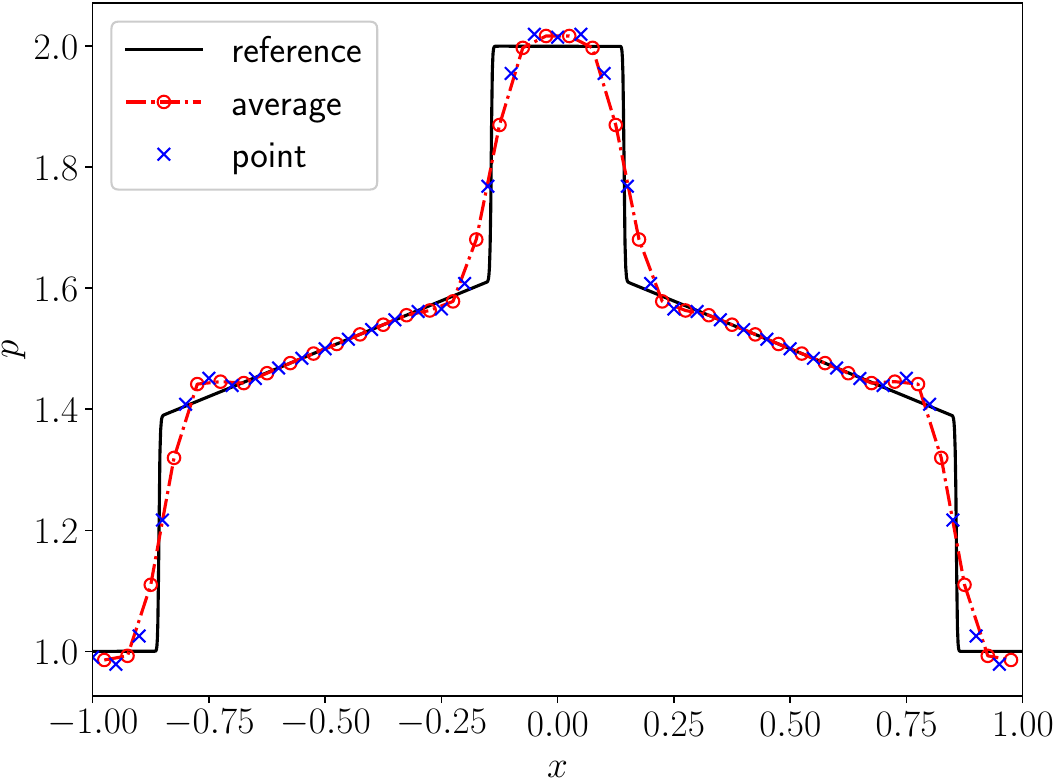}
	 	\end{subfigure}
	 	\begin{subfigure}[b]{0.465\textwidth}
	 		\centering
	 		\includegraphics[width=1.0\linewidth]{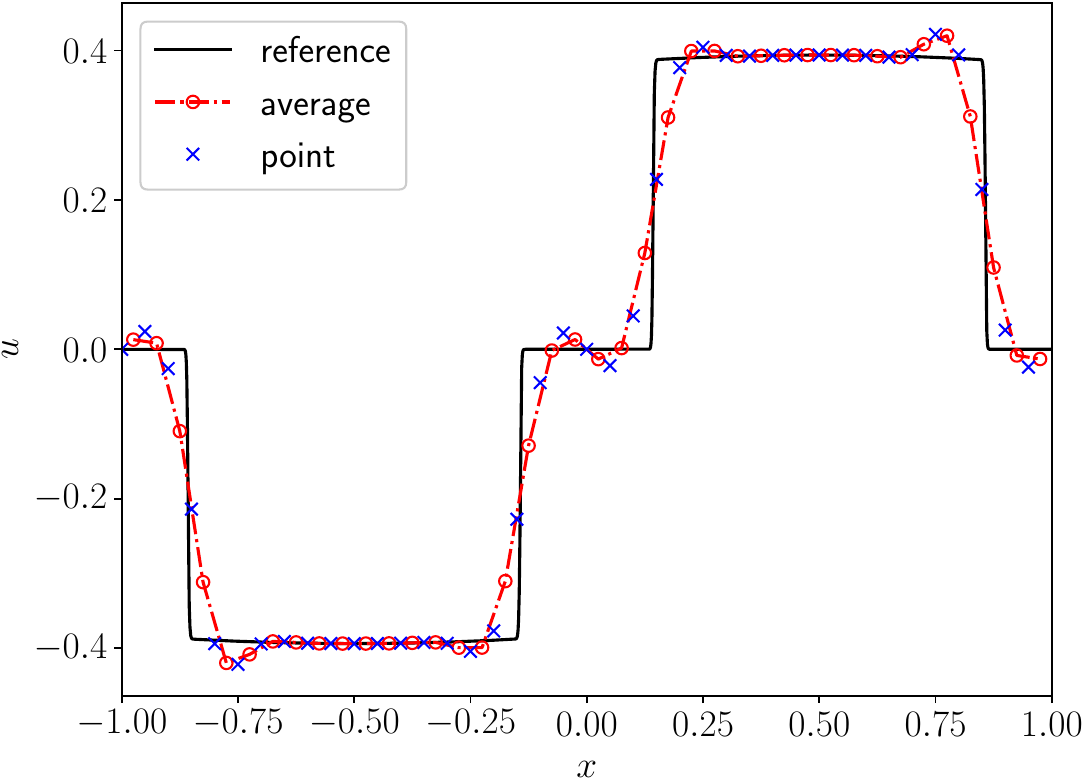}
	 	\end{subfigure}
	 	
	 	\begin{subfigure}[b]{0.45\textwidth}
	 		\centering
	 		\includegraphics[width=1.0\linewidth]{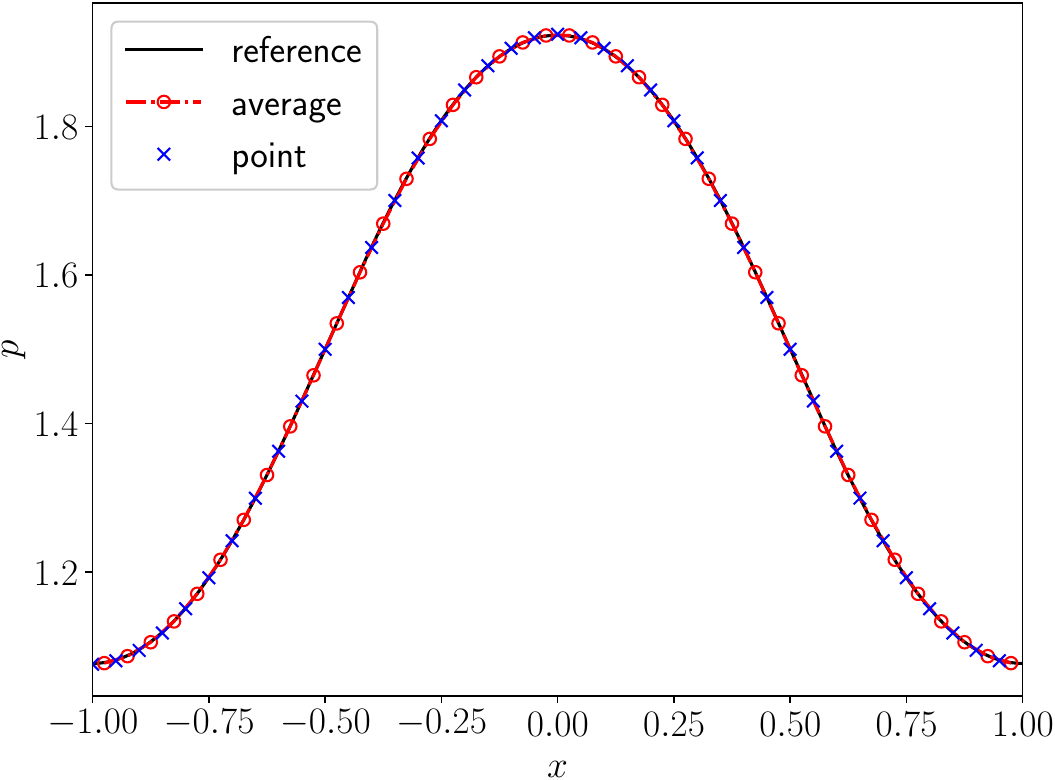}
	 	\end{subfigure}
	 	\begin{subfigure}[b]{0.465\textwidth}
	 		\centering
	 		\includegraphics[width=1.0\linewidth]{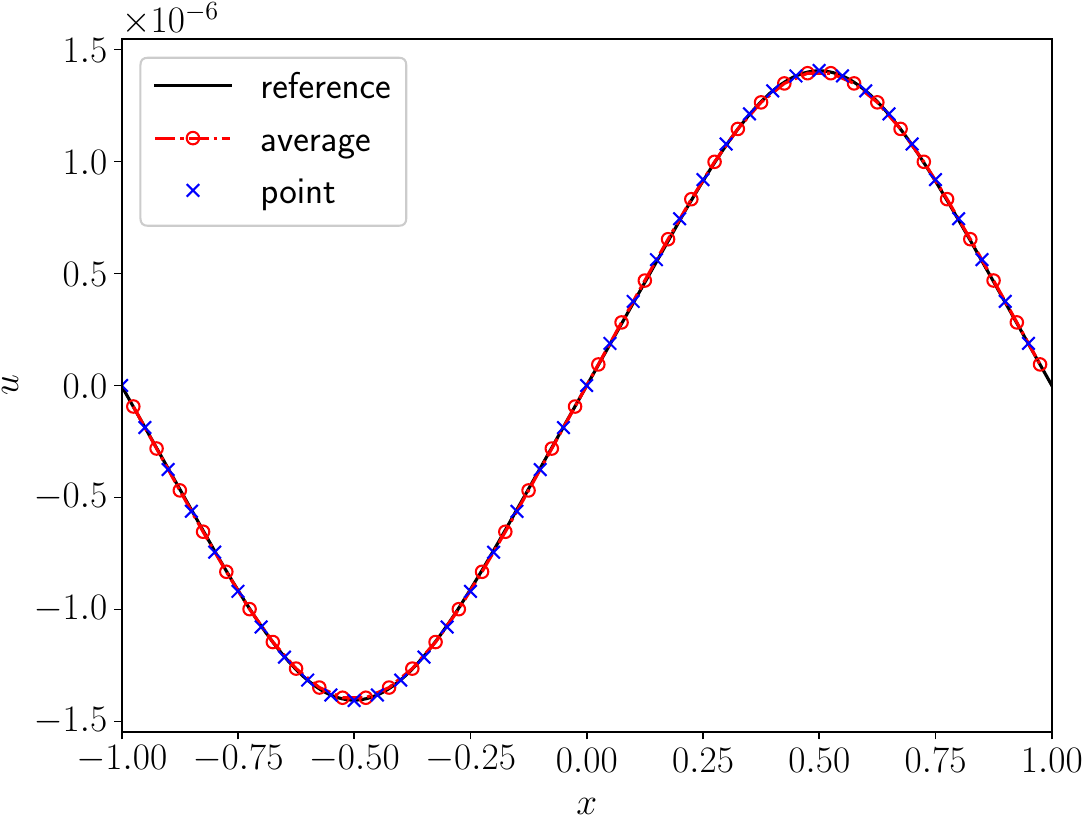}
	 	\end{subfigure}
	 	\caption{Example \ref{ex:1d_square_wave}.
	 		The results obtained by the JS-based AF scheme.
	 		Top row: $\epsilon=0.7$, bottom row: $\epsilon=10^{-6}$.
	 		Left: $p$, right: $u$.}
	 	\label{fig:1d_square_wave}
	 \end{figure}
\end{example}

\begin{example}[Variable $\sigma$]\label{ex:1d_variable_sigma}
	This test is used to show that the JS-based AF scheme also works for variable $\sigma$.
	This test problem takes the same initial condition and computational domain as Example \ref{ex:1d_square_wave},
	except for $\epsilon=1$, and $\sigma=1 + (10x)^2$.
	Note that the solution tends to be in the diffusive regime for large $\sigma$.
	
	Figure \ref{fig:1d_variable_sigma} shows the results obtained by the JS-based AF scheme on a coarse mesh of $40$ cells.
	It is seen that the scheme can capture both transport and diffusive regimes well in the domain,
	which verifies the AP property of the scheme.
	
	\begin{figure}[htbp]
		\centering
		\begin{subfigure}[b]{0.45\textwidth}
			\centering
			\includegraphics[width=1.0\linewidth]{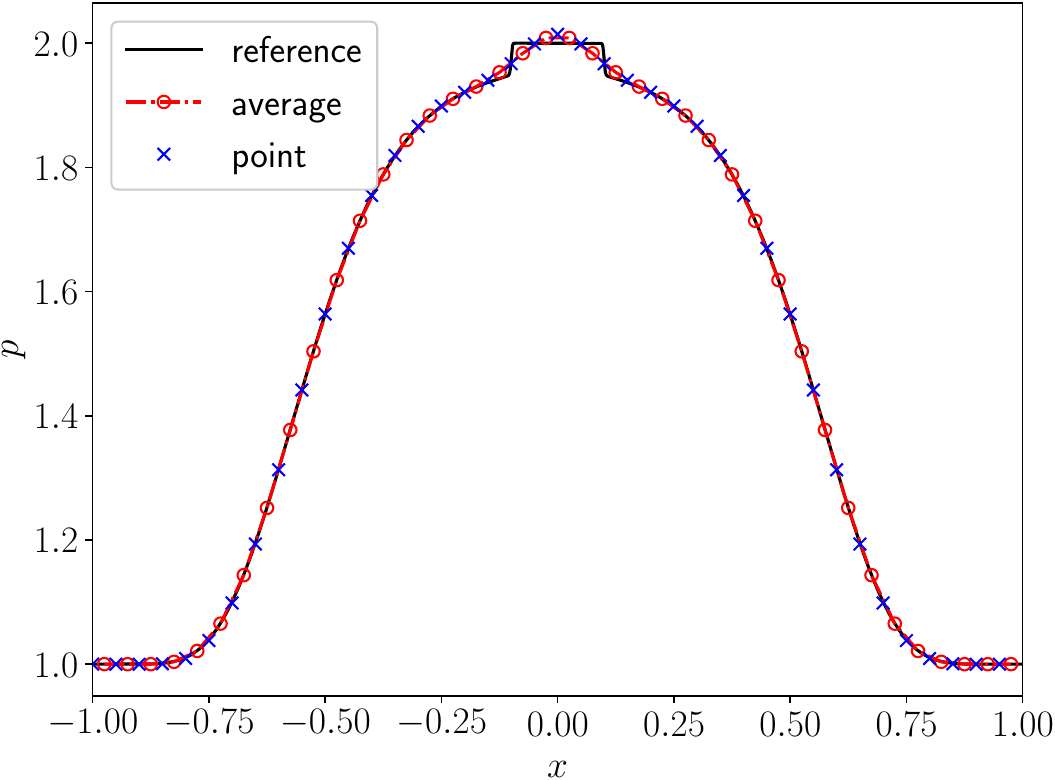}
		\end{subfigure}
		\begin{subfigure}[b]{0.465\textwidth}
			\centering
			\includegraphics[width=1.0\linewidth]{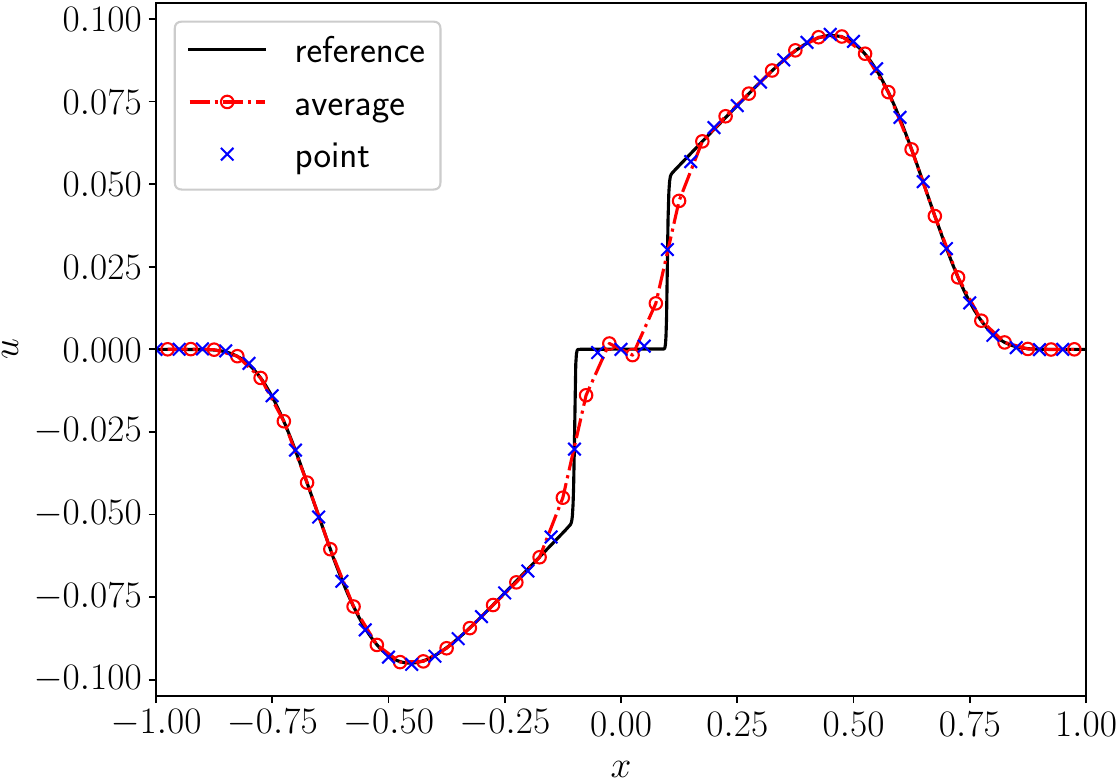}
		\end{subfigure}
		\caption{Example \ref{ex:1d_variable_sigma}.
			The results obtained by the JS-based AF scheme.
			Left: $p$, right: $u$.}
		\label{fig:1d_variable_sigma}
	\end{figure}
\end{example}
	
\begin{example}[2D accuracy test]\label{ex:2d_accuracy}
	We construct the following exact solution similar to Example \ref{ex:1d_accuracy},
	\begin{align*}
		&p = \frac{2}{r}\exp(rt)\sin(x)\sin(y),~
		u = \epsilon\exp(rt)\cos(x)\sin(y), \\
		&v = \epsilon\exp(rt)\sin(x)\cos(y),~
		r = \frac{-4}{1 + \sqrt{1-8\epsilon^2}},
	\end{align*}
	in the domain $[0,2\pi]\times [0,2\pi]$ with periodic boundary conditions and $\sigma=1$.
	The final time is $T=0.1$ and the CFL number is $0.2$ to make the spatial error dominant.	
	
	The errors and convergence orders are shown in Figure \ref{fig:2d_accuracy}.
	For $\epsilon=0.3$, the AF gets the 3rd-order accuracy for all the cell averages, also the 3rd order for the point value in $p$,
	while the 2nd order for the point value in $u,v$, which is due to the mesh alignment issue (using the LLF FVS \cite{Duan_2025_Active_SJoSC} recovers the 3rd order).	
	For $\varepsilon=10^{-6}$, all the convergence rates are 2.
	
	\begin{figure}[htbp]
		\centering
		\begin{subfigure}[b]{0.32\textwidth}
			\centering
			\includegraphics[width=1.0\linewidth]{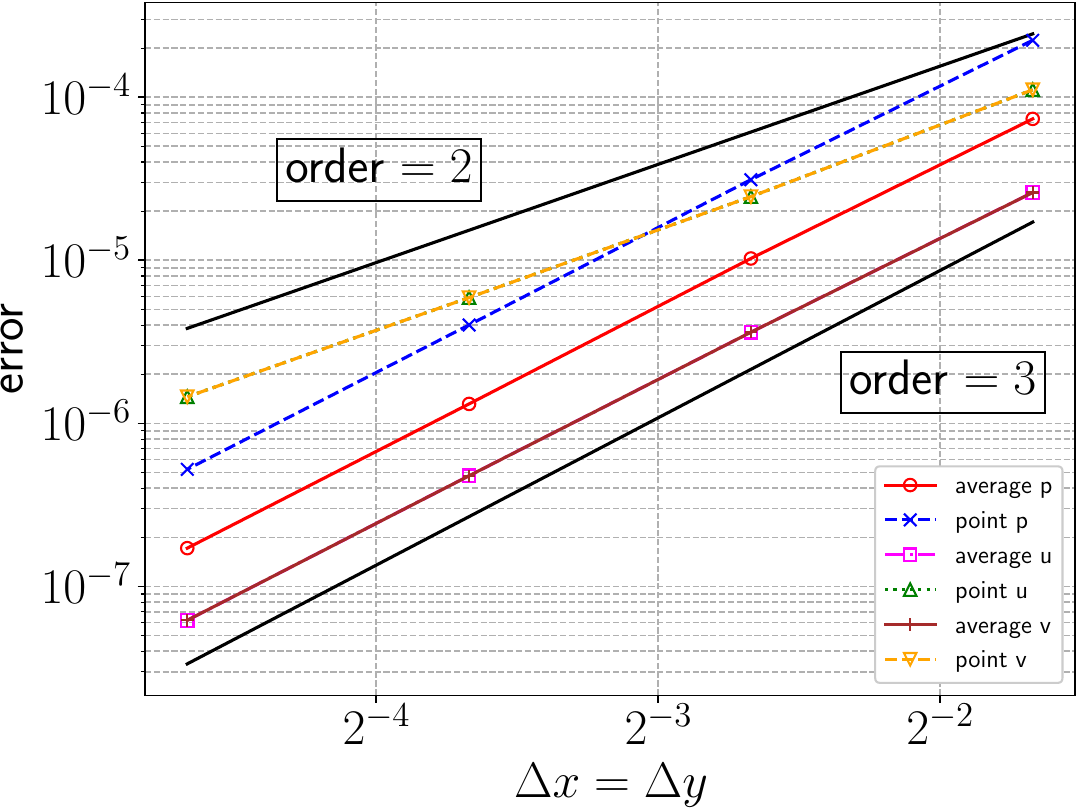}
		\end{subfigure}
		\begin{subfigure}[b]{0.32\textwidth}
			\centering
			\includegraphics[width=1.0\linewidth]{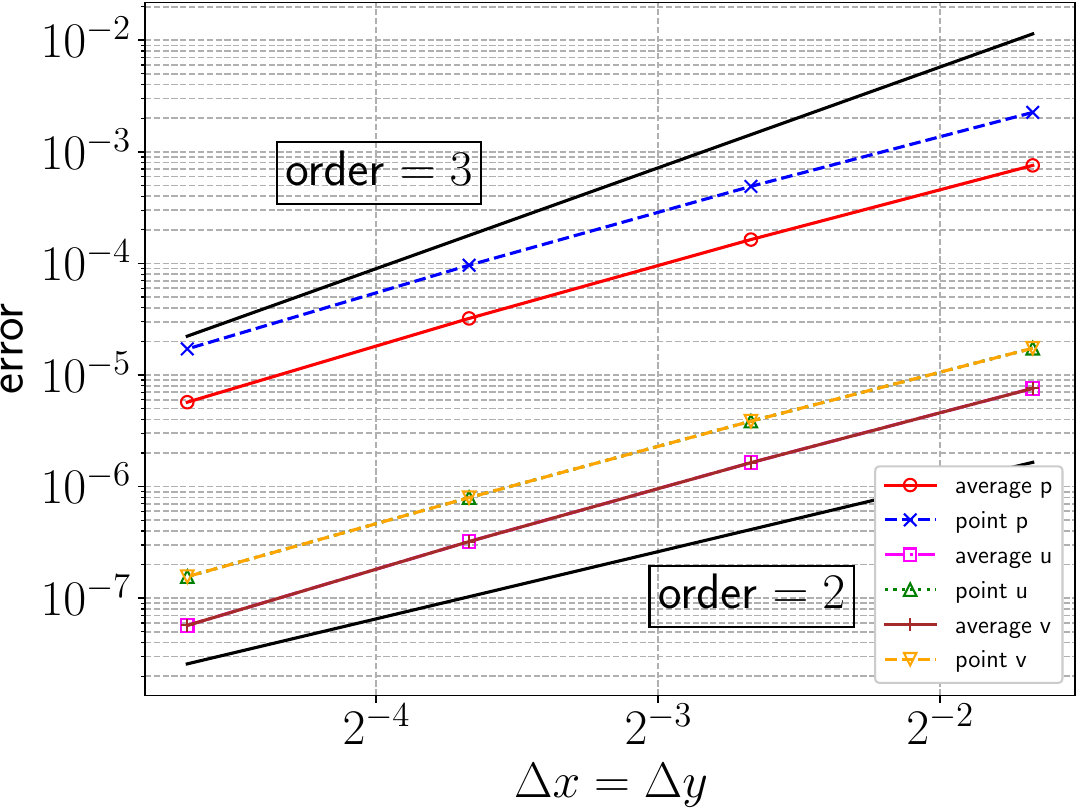}
		\end{subfigure}
		\begin{subfigure}[b]{0.325\textwidth}
			\centering
			\includegraphics[width=1.0\linewidth]{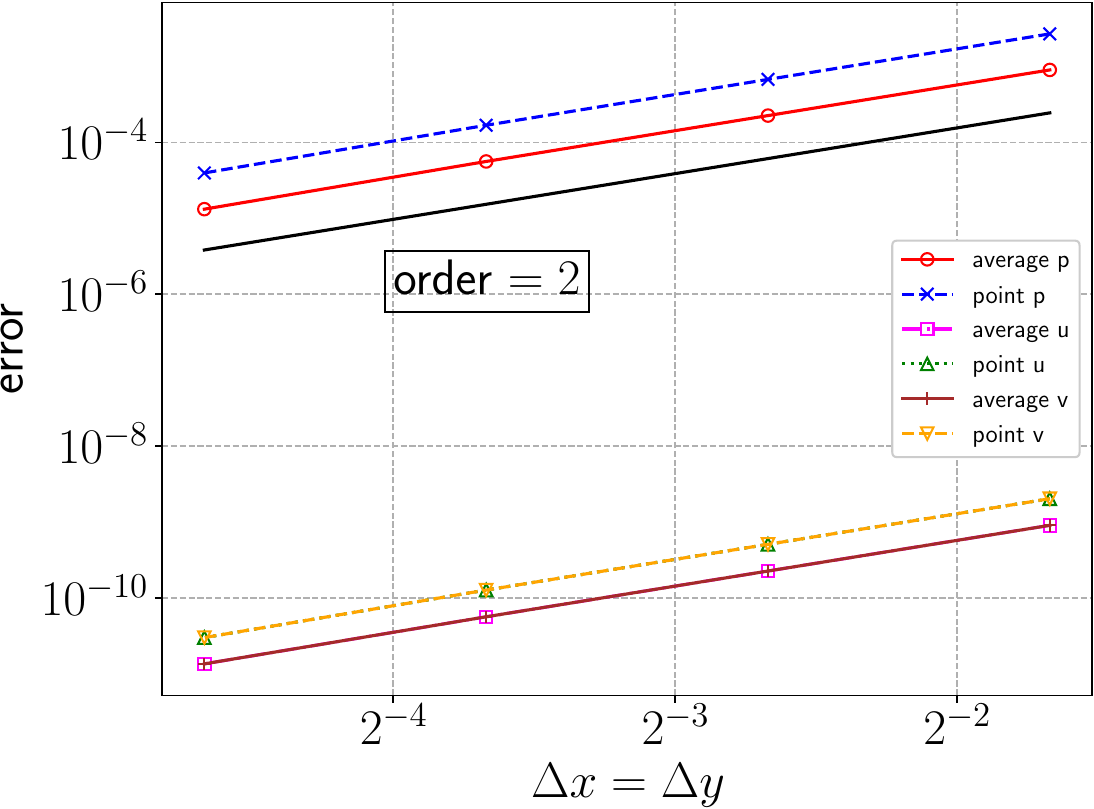}
		\end{subfigure}
		\caption{Example \ref{ex:2d_accuracy}.
			The errors and convergence rates obtained by the JS-based AF scheme and $\epsilon=0.3,10^{-2},10^{-6}$ from left to right.}
		\label{fig:2d_accuracy}
	\end{figure}
	
\end{example}

\begin{example}[Fundamental solution in transport regime]\label{ex:2d_transport_dirac}
	This test simulates wave propagation in the domain $[-1,1]\times[-1,1]$ with $\sigma=\epsilon=1$ from an initial Dirac function of $p$  \cite{Buet_2012_Design_NM},
	i.e., $p_0 = \frac{1}{\Delta x\Delta y}$ for the centered cell average otherwise $p_0 = 0$,
	and $u_0=v_0=0$ in the whole domain.
	
	Figure \ref{fig:2d_transport_dirac} shows the results obtained by using the JS-based AF scheme with $101\times 101$ cells.
	One can observe that the scheme can capture the transport phenomenon and preserve the circular shape of the wave front well.
	
	\begin{figure}[htbp]
		\centering
		\begin{subfigure}[b]{0.32\textwidth}
			\centering
			\includegraphics[width=1.0\linewidth]{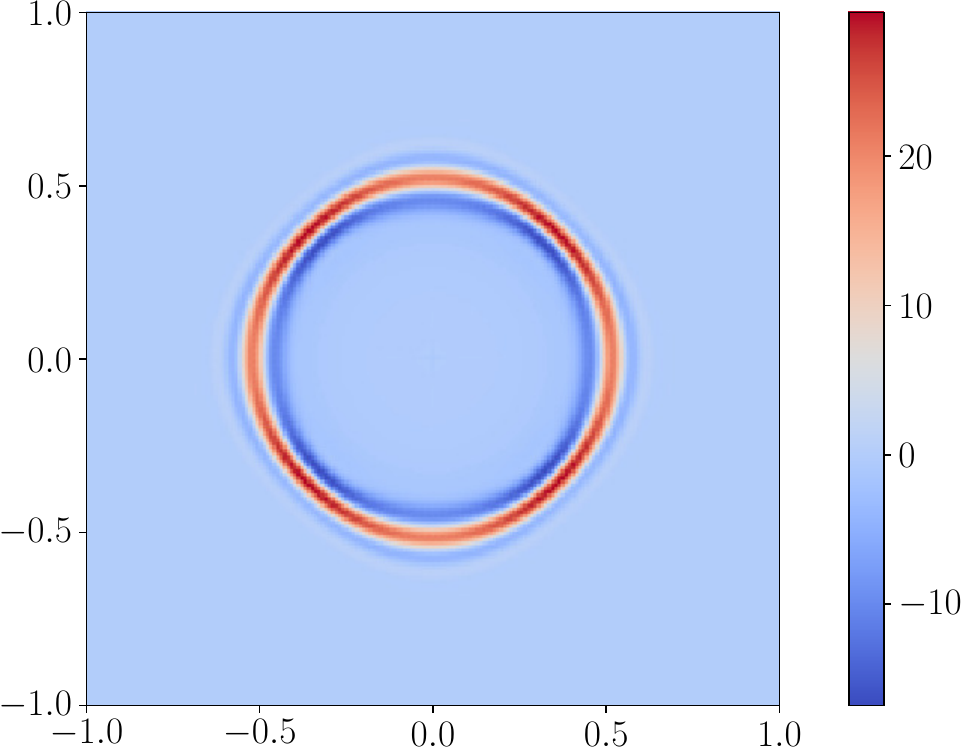}
		\end{subfigure}
		\begin{subfigure}[b]{0.32\textwidth}
			\centering
			\includegraphics[width=1.0\linewidth]{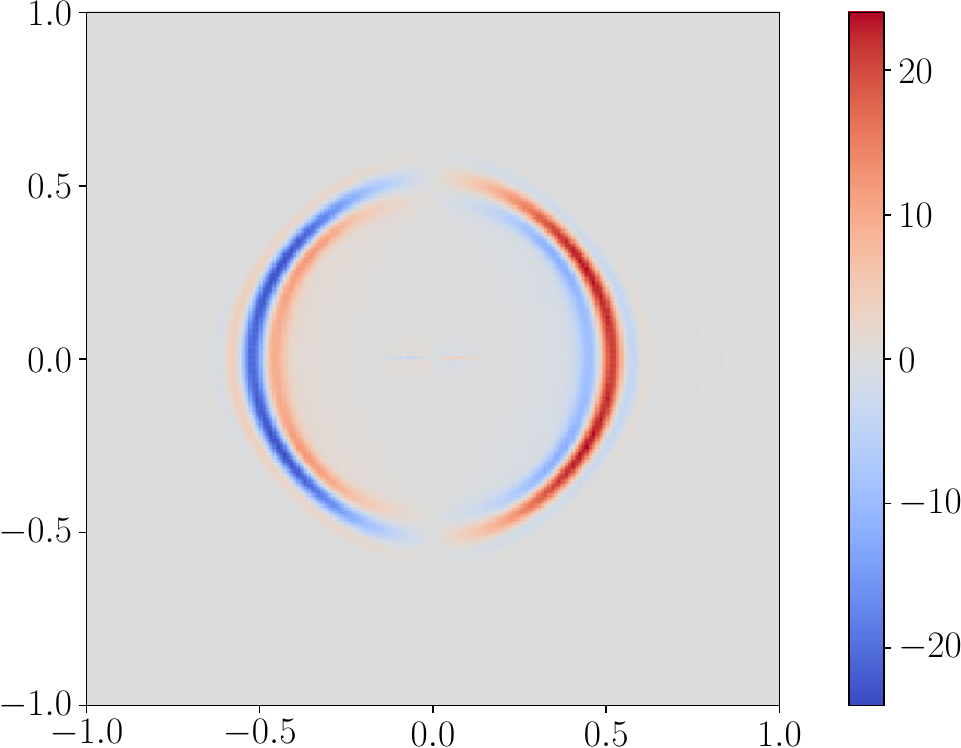}
		\end{subfigure}
		\begin{subfigure}[b]{0.32\textwidth}
			\centering
			\includegraphics[width=1.0\linewidth]{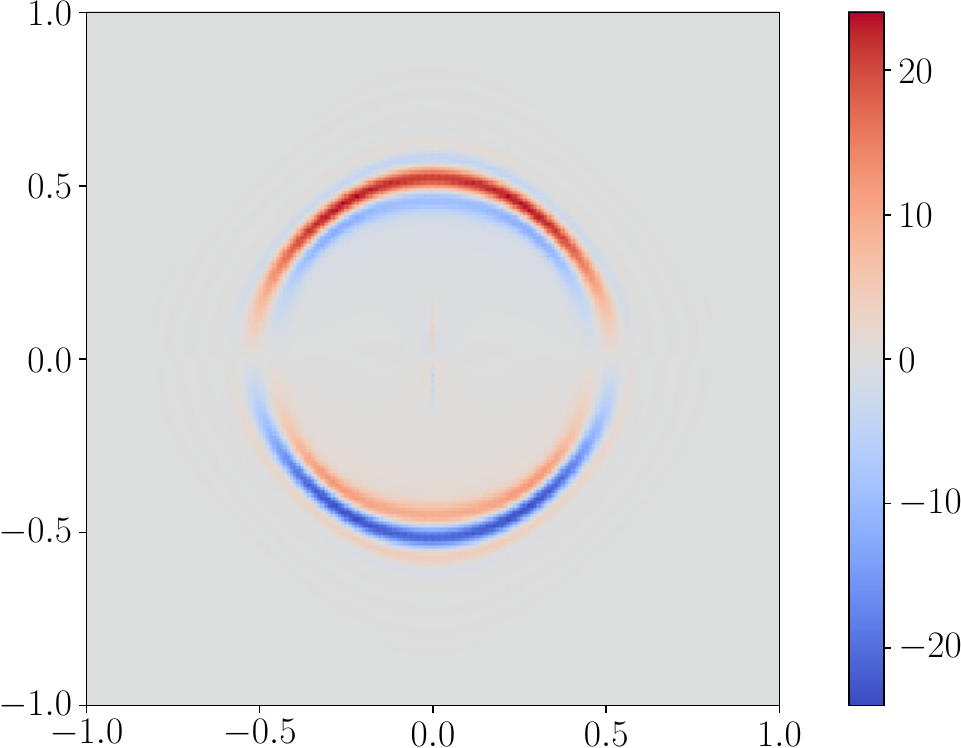}
		\end{subfigure}
		\caption{Example \ref{ex:2d_transport_dirac}.
			The results obtained by the JS-based AF scheme.
			From left to right: $p$, $u$, $v$.}
		\label{fig:2d_transport_dirac}
	\end{figure}
\end{example}

\begin{example}[Radiation]\label{ex:2d_radiation}
	This is a test problem motivated by radiation simulations,
	where $\sigma$ depends on the fluid temperature.
	We follow the setup in \cite{Buet_2012_Design_NM}.
	The domain is a square $[-1,1]\times[-1,1]$, and the initial condition is	
	\begin{equation*}
		p_0 = 10^{-3} + 100\exp(-(x^2+y^2)/0.01),~
		u_0 = v_0 = 0.
	\end{equation*}
	In this test, $\epsilon=1$, and the coefficient $\sigma$ is chosen as $1$ in the domain except for $10^4$ in eight box regions, see Figure \ref{fig:2d_radiation_sigma}.
	For example, in the first quadrant, the two boxes are $[\frac{3}{16}, \frac{7}{16}]\times [\frac{9}{16}, \frac{13}{16}]$,
	and $[\frac{9}{16}, \frac{13}{16}]\times [\frac{3}{16}, \frac{7}{16}]$,
	and the other boxes are mirror symmetric with respect to the axes.
	
	The results are obtained by using the $200\times200$ mesh, shown in Figure \ref{fig:2d_radiation}.
	The diffusion coefficient $1/\sigma$ is very small in box regions, thus the diffusion behavior is much slower, which can be clearly observed.
	Our results are comparable to those in \cite{Buet_2012_Design_NM} (using comparable mesh size),
	with better resolution near the box regions.
	
	\begin{figure}[htbp]
		\centering
		\includegraphics[width=0.3\linewidth]{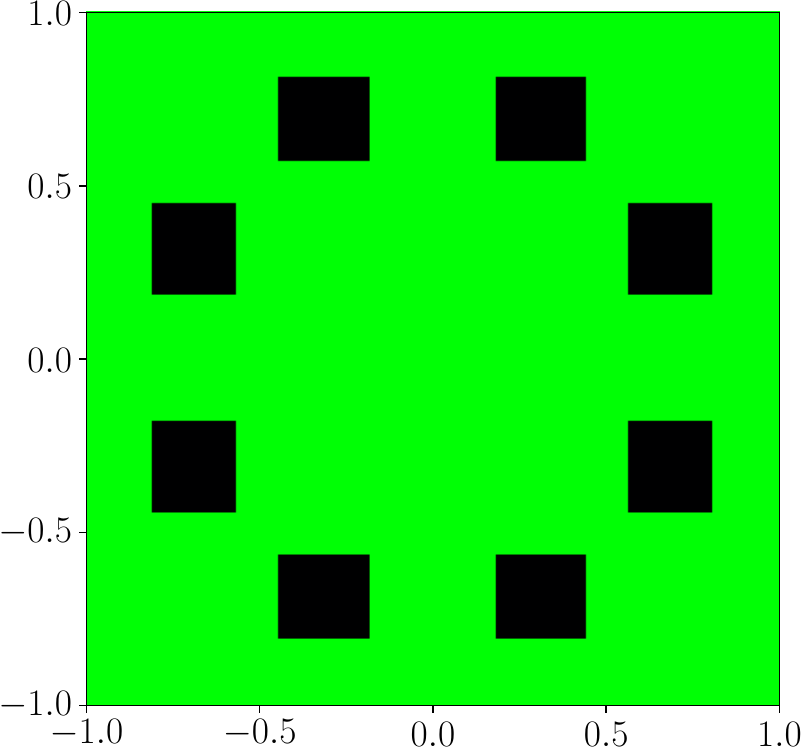}
		\caption{Example \ref{ex:2d_radiation}. $\sigma=1$ in the domain except for $\sigma=10^4$ in the black box regions.}
		\label{fig:2d_radiation_sigma}
	\end{figure}
	
	\begin{figure}[htbp]
		\centering
		\begin{subfigure}[b]{0.32\textwidth}
			\centering
			\includegraphics[width=1.0\linewidth]{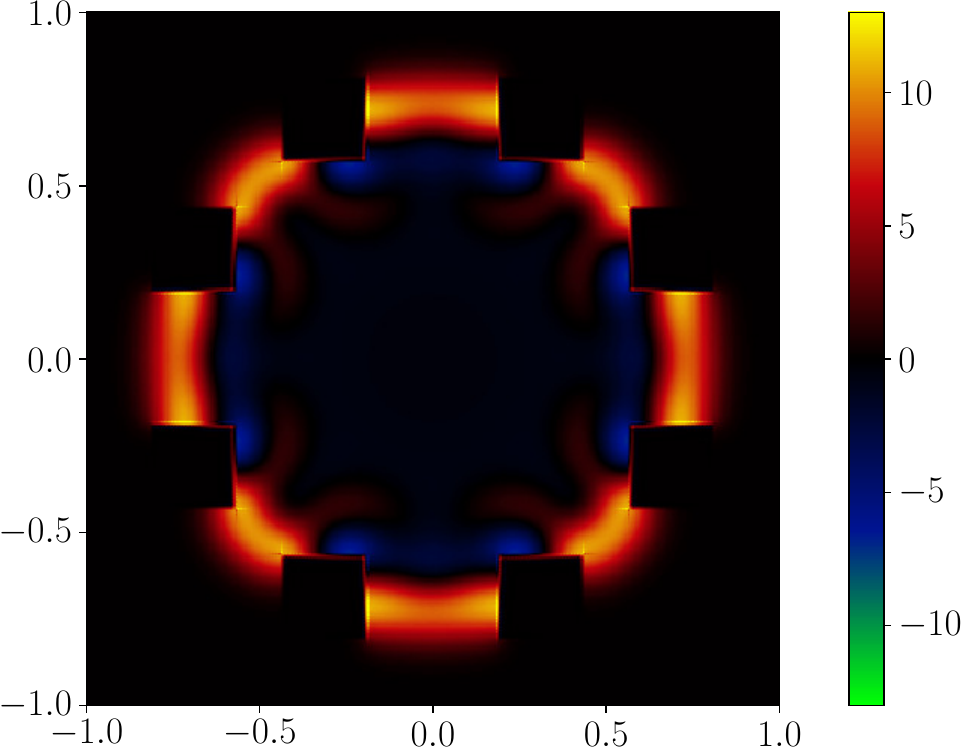}
		\end{subfigure}
		\begin{subfigure}[b]{0.32\textwidth}
			\centering
			\includegraphics[width=1.0\linewidth]{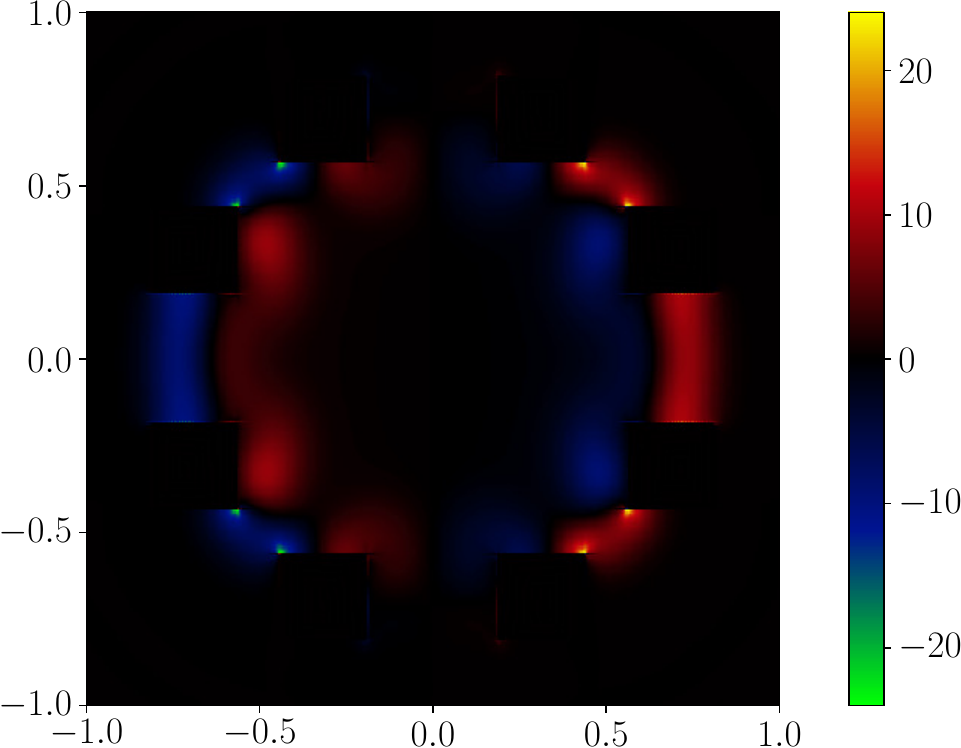}
		\end{subfigure}
		\begin{subfigure}[b]{0.32\textwidth}
			\centering
			\includegraphics[width=1.0\linewidth]{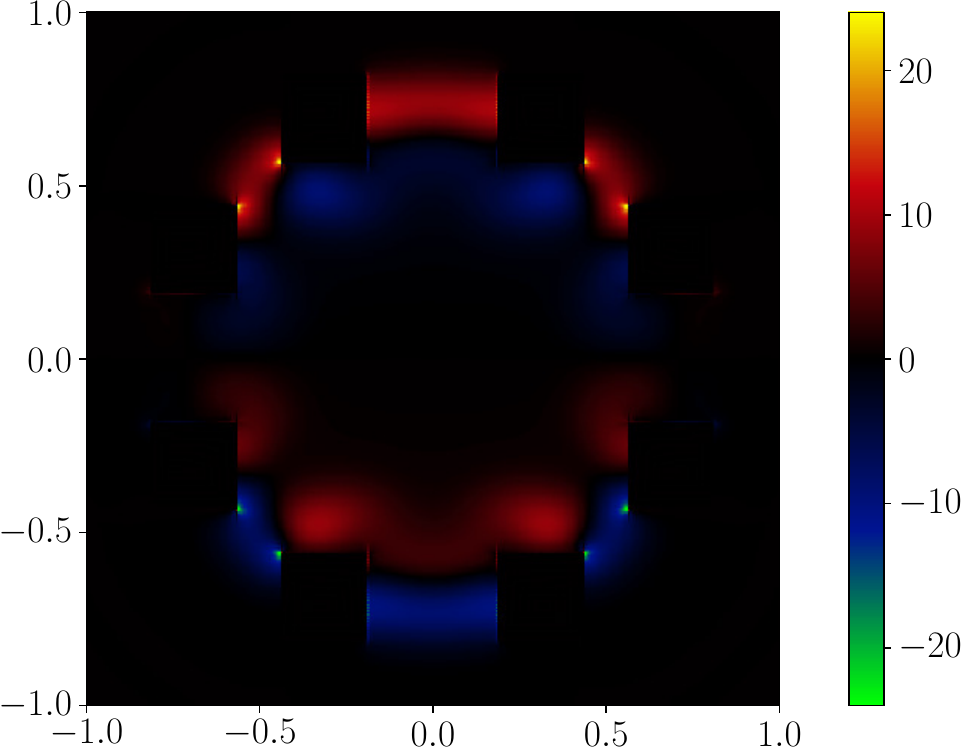}
		\end{subfigure}
		\caption{Example \ref{ex:2d_radiation}.
			The results obtained by the JS-based AF scheme.
			From left to right: $p$, $u$, $v$.}
		\label{fig:2d_radiation}
	\end{figure}
\end{example}

\section{Conclusion}\label{sec:conclusion}
In the active flux (AF) methods, the way how point values at cell interfaces are updated is essential to achieve stability, high-order accuracy, and other nice properties.
The point value update based on Jacobian splitting (JS) has been used in many existing works.
This paper investigates the JS-based AF scheme for solving the hyperbolic heat equation.
It is shown that the scheme without any modification is AP in the diffusive scaling by using the formal asymptotic analysis, discrete Fourier analysis, and numerical tests.
The convergence rates are also studied for the scheme and its limit, and the order degradation is observed in the limit.


\section*{Acknowledgement}

JD was supported by an Alexander von Humboldt Foundation Research Fellowship CHN-1234352-HFST-P. CK and WB acknowledge funding by the Deutsche Forschungsgemeinschaft (DFG, German Research Foundation) within \textit{SPP 2410 Hyperbolic Balance Laws in Fluid Mechanics: Complexity, Scales, Randomness (CoScaRa)}, project number 525941602.

\bibliographystyle{siamplain}
\bibliography{references.bib}

\newcommand{\etalchar}[1]{$^{#1}$}


\appendix

\section{Example \ref{ex:1d_accuracy} with non-well-prepared data}
\begin{figure}[htbp]
	\centering
	\begin{subfigure}[b]{0.48\textwidth}
		\centering
		\includegraphics[width=1.0\linewidth]{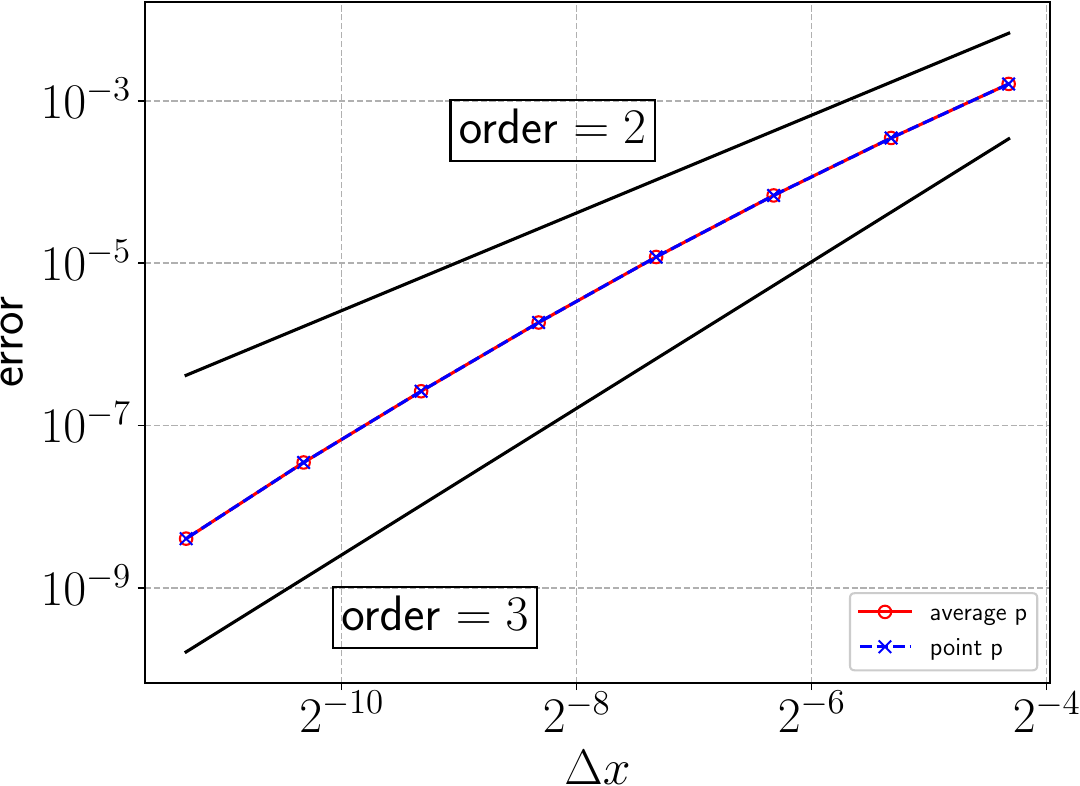}
	\end{subfigure}
	\begin{subfigure}[b]{0.48\textwidth}
		\centering
		\includegraphics[width=1.0\linewidth]{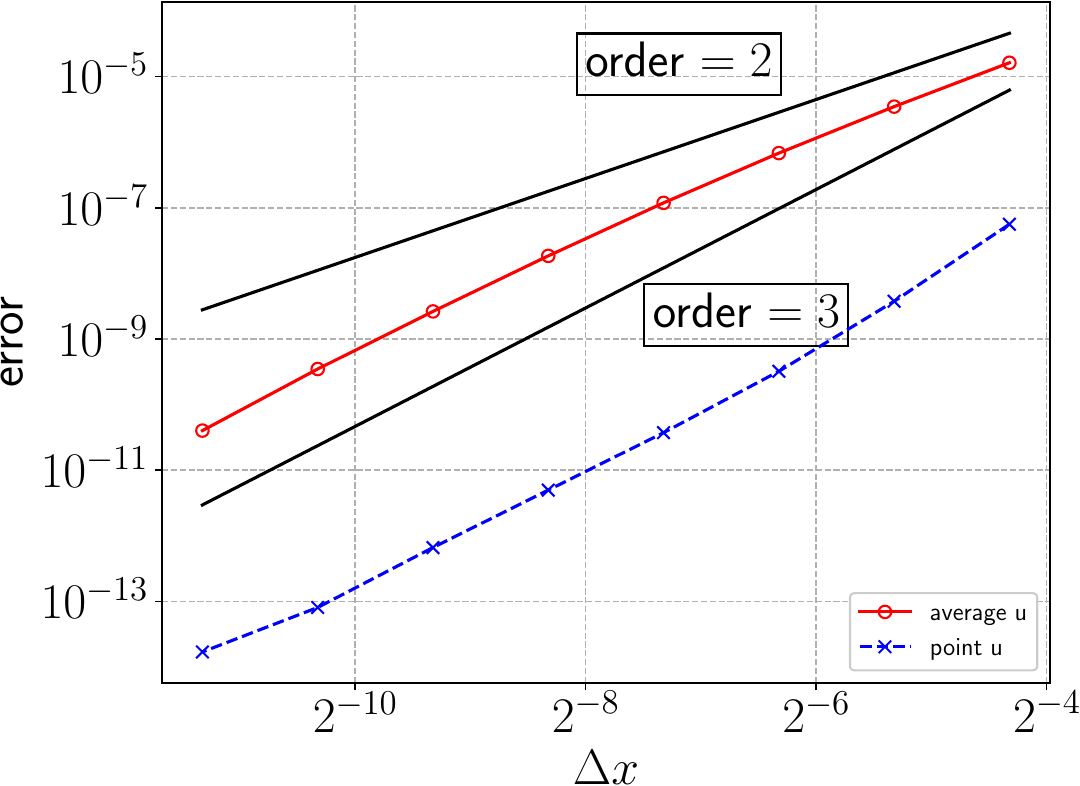}
	\end{subfigure}
	\vspace{2pt}
	
	\begin{subfigure}[b]{0.48\textwidth}
		\centering
		\includegraphics[width=1.0\linewidth]{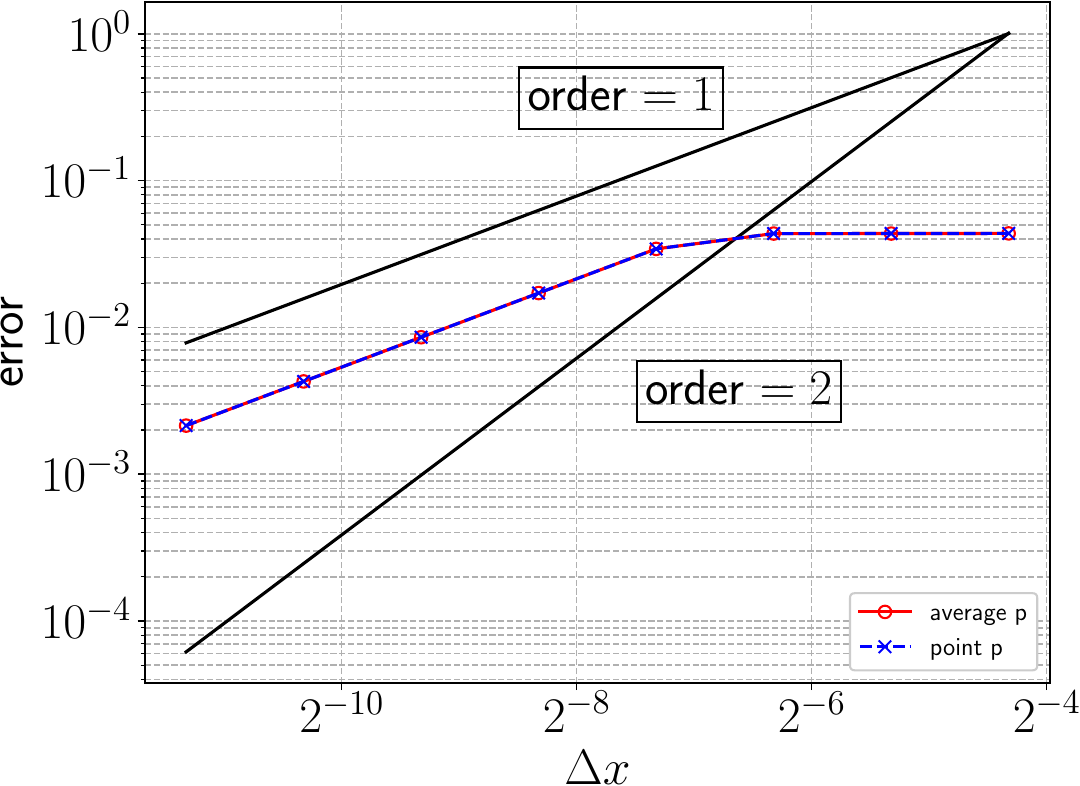}
	\end{subfigure}
	\begin{subfigure}[b]{0.48\textwidth}
		\centering
		\includegraphics[width=1.0\linewidth]{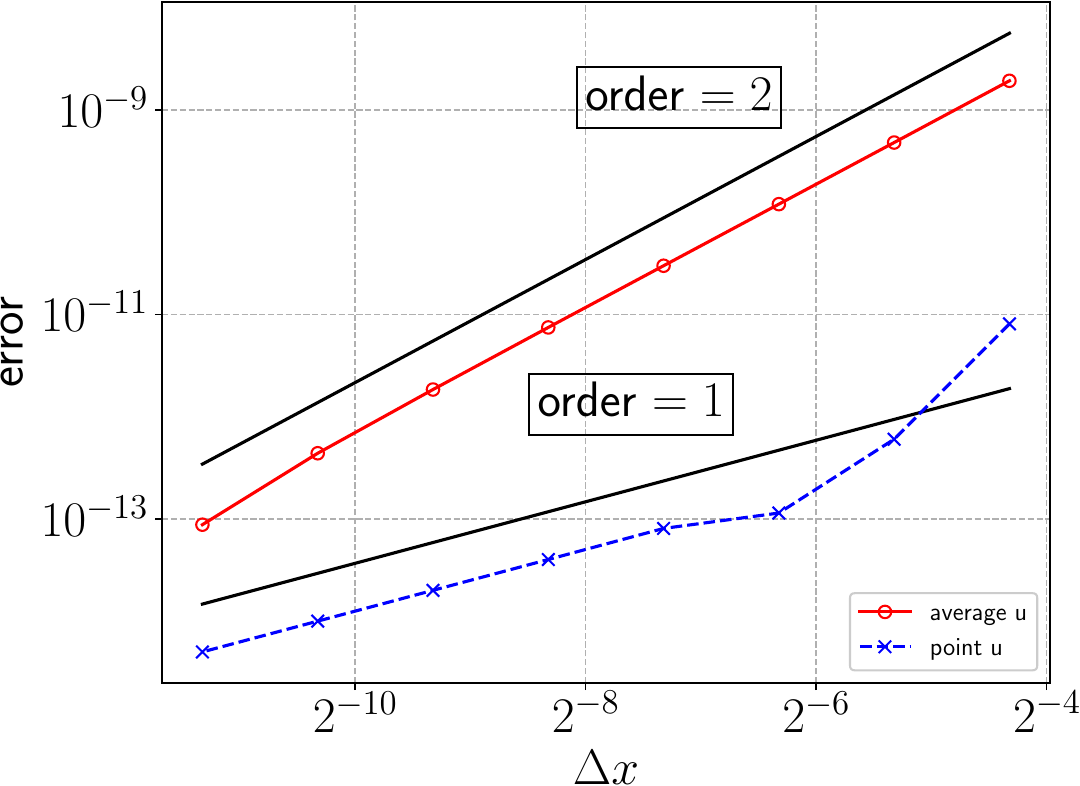}
	\end{subfigure}
	\caption{Example \ref{ex:1d_accuracy}.
		The errors and convergence rates obtained with the non-well-prepared initial data $u(0)=0.1$ and the same $p(0)$,
		where $\epsilon=10^{-2}$ (upper) and $\epsilon=10^{-6}$ (lower).}
	\label{fig:1d_accuracy_nonwellprepared}
	\end{figure}

We also compute Example \ref{ex:1d_accuracy} with non-well-prepared data $u=0.1$, and $p$ is kept the same.
The reference solution is obtained by the same scheme on a fine mesh.
Figure \ref{fig:1d_accuracy_nonwellprepared} shows the errors and convergence rates. 
The 3rd-order convergence is observed for $\epsilon=10^{-2}$ when $\Delta x$ is small enough,
while for $\epsilon=10^{-6}$, the convergence rate is the 1st order except for the 2nd-order accuracy in the cell average of $u$.

\end{document}